\documentclass[reqno,11pt]{amsart}
\usepackage{amsmath,mathrsfs}
\usepackage{amssymb, amsmath}
\usepackage{graphicx}
\usepackage{pdfsync}
\usepackage{color}
\usepackage[colorlinks,
linkcolor=red,
anchorcolor=blue,
citecolor=green
]{hyperref}
\usepackage{colonequals}
\def\inte#1{
	\displaystyle\mathop{#1\kern0pt}^\circ }

\let\f=\frac

\def\virgp{\raise 2pt\hbox{,}}
\def\cdotpv{\raise 2pt\hbox{;}}

\def\C{\mathop{\mathbb C\kern 0pt}\nolimits}
\def\DD{\mathop{\mathbb D\kern 0pt}\nolimits}
\def\EE{\mathop{{\mathbb E \kern 0pt}}\nolimits}
\def\K{\mathop{\mathbb K\kern 0pt}\nolimits}
\def\N{\mathop{\mathbb N\kern 0pt}\nolimits}
\def\Q{\mathop{\mathbb Q\kern 0pt}\nolimits}
\def\R{\mathop{\mathbb R\kern 0pt}\nolimits}
\def\SS{\mathop{\mathbb S\kern 0pt}\nolimits}
\def\ZZ{\mathop{\mathbb Z\kern 0pt}\nolimits}
\def\TT{\mathop{\mathbb T\kern 0pt}\nolimits}
\def\P{\mathop{\mathbb P\kern 0pt}\nolimits}

\newcommand{\beq}{\begin{equation}}
	\newcommand{\eeq}{\end{equation}}
\newcommand{\ben}{\begin{eqnarray}}
	\newcommand{\een}{\end{eqnarray}}
\newcommand{\beno}{\begin{eqnarray*}}
	\newcommand{\eeno}{\end{eqnarray*}}

\newtheorem{thm}{Theorem}[section]
\newtheorem{lem}{Lemma}[section]

\theoremstyle{definition}

\newtheorem{rmk}{Remark}[section]

\numberwithin{equation}{section}

\begin{document}
	\title[Coercivity estimate]
	{An explicit coercivity estimate of the linearized quantum Boltzmann operator without angular cutoff}

	\author[T. Yang]{Tong Yang}
	\address[T. Yang]{Department of Applied Mathematics, 
		The Hong Kong Polytechnic University,
		Hong Kong, P.R. China.}
	\email{t.yang@polyu.edu.hk}
	
	\author[Y.-L. Zhou]{Yu-long Zhou}
	\address[Y.-L. Zhou]{School of Mathematics, Sun Yat-Sen University, Guangzhou, 510275, P. R.  China.} \email{zhouyulong@mail.sysu.edu.cn}

	\begin{abstract} The quantum Boltzmann-Bose equation describes a large system of Bose-Einstein
		particles in the  weak-coupling regime. If the particle interaction is governed
		by the inverse power law, the corresponding collision  kernel has angular singularity. 
		In this paper, we give a constructive proof of the coercivity estimate for the linearized quantum Boltzmann-Bose operator to capture the effects of the singularity and the fugacity. 
		Precisely, the estimate explicitly reveals the dependence on  the
		 fugacity parameter  before the Bose-Einstein condensation. 
		With the coercivity estimate, the global in time well-posedness of the inhomogeneous quantum Boltzmann-Bose equation in the
		perturbative framework and stability of the Bose-Einstein equilibrium can be established. 
	\end{abstract}

	\maketitle
	\markright{Quantum Boltzmann equation}

	\setcounter{tocdepth}{1}
	

	
	
	\noindent {\sl\small AMS Subject Classification (2020):} {35Q20, 82C40.}


	%

	\section{Introduction}
	
	%
	%
	%
	%
	
	Quantum Boltzmann equations are proposed to describe the time evolution of a large system of weakly interacting bosons or fermions. The derivation of such equations dates back to  1928 by Nordheim \cite{nordhiem1928kinetic} and 1933 by Uehling-Uhlenbeck \cite{uehling1933transport}.
	As a result, the quantum Boltzmann equations are also called Boltzmann-Nordheim equations or Uehling-Uhlenbeck equations.
	Later on, further developments were made by Erd{\H{o}}s-Salmhofer-Yau \cite{erdHos2004quantum}, Benedetto-Castella-Esposito-Pulvirenti \cite{benedetto2004some}, and \cite{benedetto2006some,benedetto2007short,benedetto2008n,benedetto2005weak}, also Lukkarinen-Spohn \cite{lukkarinen2009not}. See \cite{chen2021emergence} for latest derivation with error analysis.
	One can refer the classical book \cite{chapman1990mathematical} for 
	description on its physical background.
	
	The inhomogeneous quantum Boltzmann equation reads
	\ben \label{quantum-Boltzmann-UU}
	\partial _t F +  v \cdot \nabla_{x} F=Q_{\Phi, \hbar}(F,F), ~~t > 0, x \in \Omega, v \in \R^3.
	\een
	Here $F(t,x,v)\geq 0$ is the density
	function of particles with velocity
	$v\in\R^3$ at time $t\geq 0$ in position $x \in \Omega$. 
	To avoid the difficulty from the boundary, one can study the Cauchy problem with initial data $F|_{t=0}(x,v) = F_{0}(x,v)$
	 in either a  periodic box $\Omega = \mathbb{T}^{3} \colonequals  [0,1]^{3}$ or in the whole space $\Omega = \mathbb{R}^{3}$. 
	
	The  quantum Boltzmann operator $Q_{\Phi, \hbar}$ for bosons is given by
	\ben  \nonumber
	&& Q_{\Phi, \hbar}(g,h)(v) 
	\\ \label{U-U-operator} &\colonequals&   \int_{{\mathbb S}^2 \times {\mathbb R}^3} B_{\Phi,\hbar}(v- v_{*},\sigma)
	\big(g_{*}^{\prime} h^{\prime}(1 +  \hbar^{3}g_{*})(1 +  \hbar^{3}h) - g_{*} h(1 +  \hbar^{3}g_{*}^{\prime})(1 +  \hbar^{3}h^{\prime})\big)
	\mathrm{d}\sigma \mathrm{d}v_{*},
	\een
	where according to \cite{erdHos2004quantum} and \cite{benedetto2005weak} the Boltzmann kernel $B_{\Phi,\hbar}(v- v_{*},\sigma)$ has the following form
	\ben \label{scaling-Boltzmann-kernel}
	B_{\Phi,\hbar}(v- v_{*},\sigma) \colonequals   \hbar^{-4} |v-v_{*}| \big(
	\hat{\Phi} (\hbar^{-1} |v-v^{\prime}|)
	+ \hat{\Phi} (\hbar^{-1} |v-v^{\prime}_{*}|)
	\big)^{2}.
	\een
	Here $\hbar$ is the Plank constant and $\Phi$ is a radial potential function. The Fourier transform $\hat{\Phi}$ of $\Phi$ is still a radial function so that it depends only on the magnitude of its argument in \eqref{scaling-Boltzmann-kernel}. 	Here,
	we use the conventional notations $h=h(v)$, $g_*=g(v_*)$,
	$h'=h(v')$, $g'_*=g(v'_*)$ where $v'$, $v_*'$ are given by
	\ben\label{v-prime-v-prime-star}
	v'=\frac{v+v_{*}}{2}+\frac{|v-v_{*}|}{2}\sigma, \quad  v'_{*}=\frac{v+v_{*}}{2}-\frac{|v-v_{*}|}{2}\sigma, \quad \sigma\in\SS^{2}.
	\een
	
%

	%

	By using the following scaling
	\ben \label{scaling-tranform}
	\tilde{F}(t,x,v) = \hbar^{3}F(\hbar^{3}t,x,\hbar^{-3}v), \quad \phi(|x|) = \hbar^{4}\Phi(\hbar^{4}|x|),
	\een
	one can normalize the Plank constant $\hbar$. Indeed, it is easy to check $F$ is a solution to \eqref{quantum-Boltzmann-UU} if and only if $\tilde{F}$ is a solution of the following normalized equation
	\ben \label{quantum-Boltzmann-UU-scaling}
	\partial _t F +  v \cdot \nabla_{x} F=Q_{\phi}(F,F), ~~t > 0, x \in \Omega, v \in \R^3,
	\een
	where the operator $Q_{\phi}$ is defined by
	\ben \label{U-U-operator-scaling}
	Q_{\phi}(g,h)(v) \colonequals   \int_{{\mathbb S}^2 \times {\mathbb R}^3} B_{\phi}(v- v_{*},\sigma)
	\big(g_{*}^{\prime} h^{\prime}(1 +  g_{*})(1 +  h) - g_{*} h(1 +  g_{*}^{\prime})(1 +  h^{\prime})\big)
	\mathrm{d}\sigma \mathrm{d}v_{*}
	\een
	with the kernel $B_{\phi}(v- v_{*},\sigma)$ given by
	\ben \label{scaling-Boltzmann-kernel-scaling}
	B_{\phi}(v- v_{*},\sigma) \colonequals   |v-v_{*}| \big(
	\hat{\phi} (|v-v^{\prime}|)
	+ \hat{\phi} (|v-v^{\prime}_{*}|)
	\big)^{2}.
	\een
	
	\subsection{Potential function and Boltzmann kernel} The mathematical study of the BBE (Boltzmann-Bose-Einstein)
	 equation \eqref{quantum-Boltzmann-UU-scaling} is  difficult due to the high order nonlinearity and possible blow-up in low temperature
	that corresponds to the intriguing  Bose-Einstein condensation (BEC).  As a result, most of the existing results on the BBE equation are concerned with hard sphere model $B_{\phi} = (v- v_{*},\sigma) = C |v-v_{*}|$ for $\phi(x) = \delta(x)$ in \eqref{scaling-Boltzmann-kernel-scaling}. In the homogeneous case, one can refer to \cite{lu2000modified,lu2004isotropic,escobedo2007fundamental,escobedo2008singular} for existence of
	isotropic solutions $F(t,v)=F(t,|v|)$. As for anisotropic solutions, the work \cite{briant2016cauchy} is for  local existence and \cite{li2019global} for global existence. In the inhomogeneous case,
	we refer to the two recent results \cite{bae2021relativistic,ouyang2022quantum} for global well-posedness in the perturbative framework.

	Another physically important potential is the inverse power law $\phi(x) = |x|^{-p}$ which has been extensively studied for the classical Boltzmann equation but much  less is known for the quantum Boltzmann equation. 
	Note that $p=1$ corresponds to the  Coulomb potential that is the critical value for Boltzmann equation to be well-defined.
	In this paper, we work in three dimensional space so that   $p=3$ is the critical value for $|x|^{-p}$
	 to be locally integrable near $|x|=0$. Fix $1<p<3$. Then the Fourier transform of $\phi(x) = |x|^{-p}$ is
	$
	\hat{\phi}(\xi) = \hat{\phi}(|\xi|) = C |\xi|^{p-3}.
	$ Let $\theta$ be the angle between $v-v_{*}$ and $\sigma$, then
	$ |v-v^{\prime}|=|v-v_{*}| \sin \frac{\theta}{2},  |v-v^{\prime}_{*}|=|v-v_{*}| \cos \frac{\theta}{2} $. By \eqref{scaling-Boltzmann-kernel-scaling}, the kernel then reads
	\beno
	B_{\phi}(v- v_{*},\sigma) = C |v-v_{*}|^{2p-5} (\sin^{p-3} \frac{\theta}{2} + \cos^{p-3} \frac{\theta}{2})^{2}.
	\eeno
	Throughout the paper,  $C$ denotes a  generic constant.
	Thanks to the symmetry of \eqref{scaling-Boltzmann-kernel-scaling}, one can  assume $0 \leq \theta \leq \pi/2$. Then the fact that $\sin \frac{\theta}{2} \leq \cos \frac{\theta}{2}$ and $p<3$ 
	implies that 
	\beno
	B_{\phi}(v- v_{*},\sigma) \leq C |v-v_{*}|^{2p-5} \sin^{2p-6} \frac{\theta}{2}.
	\eeno
	Hence, for $p>1$,  
	\ben \label{ub-ipl-mean-momentum-transfer}
	\int B_{\phi}(v- v_{*},\sigma) \sin^{2} \frac{\theta}{2} \mathrm{d}\sigma \leq C |v-v_{*}|^{2p-5}.
	\een
	Note that the integral domain $\sigma \in \mathbb{S}^{2}$ is omitted for brevity.  When there is no ambiguity, the range of some frequently
	used variables will be omitted from now on. To be  consistent, we always take
	 $\sigma \in \mathbb{S}^{2}, x \in \Omega, v, v_{*} \in \mathbb{R}^{3}$ unless otherwise specified.

	By \eqref{ub-ipl-mean-momentum-transfer}, for $1<p<3$ the kernel $B_{\phi}$ has finite mean momentum transfer which is a necessary condition for the classical Boltzmann operator to be well-defined.
	The constant $C$ in \eqref{ub-ipl-mean-momentum-transfer} blows up as $p \rightarrow 1^{+}$ since $\int_{0}^{\sqrt{2}/{2}} t^{2p-3} \mathrm{d}t \sim \frac{1}{p-1}$ for $1<p<3$.
	If $p > 2$, the angular function $\sin^{2p-6} \frac{\theta}{2}$ satisfies Grad's angular cutoff assumption
	\beno
	\int \sin^{2p-6} \frac{\theta}{2} \mathrm{d}\sigma \lesssim  \frac{1}{p-2}.
	\eeno
	Note that $p=2$ is the critical value such that Grad's cutoff assumption fails. Hence,
	$2<p<3$ can be viewed as angular cutoff while $1<p<2$ as angular non-cutoff.
	We consider the more difficult case when $1<p<2$ in this paper. For this, by taking $s=2-p, \gamma=2p-5$, it suffices to consider the following more general kernel
	\ben \label{Boltzmann-kernel-general-and-old}
	B_{\phi} (v- v_{*},\sigma)  =  C |v-v_{*}|^{\gamma} (\sin^{-1-s} \frac{\theta}{2} + \cos^{-1-s} \frac{\theta}{2})^{2} \mathrm{1}_{0 \leq \theta \leq \pi/2}, \quad -3<\gamma \leq  0<s<1.
	\een
	The parameter pair $(\gamma, s)$ is commonly used in the study of classical Boltzmann equations with inverse power law. Indeed,
	motivated by the inverse power law potential, the Boltzmann kernel can be written as 
	\ben \label{inverse-power-law-kernel}
	B(v-v_{*},\sigma) = B^{\gamma,s}(v-v_{*},\sigma) \colonequals    |v-v_{*}|^{\gamma} \sin^{-2-2s}\frac{\theta}{2} \mathrm{1}_{0 \leq \theta \leq \pi/2}.
	\een
	Since $\sin \frac{\theta}{2} \leq \cos \frac{\theta}{2}$ for $0 \leq \theta \leq \frac{\pi}{2}$, the two kernel functions \eqref{inverse-power-law-kernel} and \eqref{Boltzmann-kernel-general-and-old} satisfy
	\ben \label{Boltzmann-kernel-general-ub}
	C B^{\gamma,s}(v-v_{*},\sigma) \leq B_{\phi}(v- v_{*},\sigma) \leq 4 C B^{\gamma,s}(v-v_{*},\sigma) .
	\een
	Hence, we will focus on the kernel  function given in \eqref{inverse-power-law-kernel}. 
	From now on, the notation $B(v- v_{*},\sigma)$ or $B$ stands for the kernel given in \eqref{inverse-power-law-kernel} unless otherwise specified. 	
	The condition $\gamma+2s < 0$ is referred to as soft potentials, and $\gamma+2s\ge 0$ as hard potentials.
	
	\subsection{Classical Boltzmann equation}
	The well-posedness of classical Boltzmann  was established in  \cite{alexandre2011boltzmann,alexandre2012boltzmann} in the whole space  and  \cite{gressman2011global} in torus  independently for the non-cutoff kernel $B^{\gamma,s}$ \eqref{inverse-power-law-kernel} in the perturbative framework. Let us briefly recall the main estimates as follows. The classical Boltzmann equation reads 
\ben \label{classical-Boltzmann-UU}
\partial _t F +  v \cdot \nabla_{x} F = Q^{\gamma,s}_{c}(F,F),
\een	
	where 
	\ben \label{classical-Boltzmann-operator}
	Q^{\gamma,s}_{c}(g,h)(v) \colonequals \int B^{\gamma,s}(v- v_{*},\sigma)
	\big(g_{*}^{\prime} h^{\prime} - g_{*} h \big) \mathrm{d}\sigma \mathrm{d}v_{*}
	\een
	The subscript
	``$c$'' here  refers to ``classical''. Recall that the global Maxwellian density
	 $\mu (v) \colonequals  e^{-\f{1}{2}|v|^{2}}$ is a typical equilibrium of  \eqref{classical-Boltzmann-UU}. 
	By setting $F = \mu + \mu^{\f12}f$, the equation \eqref{classical-Boltzmann-UU} is written as
	\beno 
	\partial _t f +  v \cdot \nabla_{x} f + \mathcal{L}^{\gamma,s}_{c}f = \Gamma^{\gamma,s}_{c}(f,f)
	\eeno
	where the linearized Boltzmann operator $\mathcal{L}^{\gamma,s}_{c}$ and the nonlinear term $\Gamma^{\gamma,s}_{c}$ are given by
		\ben
	\label{Def-L-and-Gamma}
 \mathcal{L}^{\gamma,s}_{c}f := - \mu^{-\f12} Q^{\gamma,s}_{c}(\mu,\mu^{\f12}f) - \mu^{-\f12} Q^{\gamma,s}_{c}(\mu^{\f12}f, \mu), \quad 
	 \Gamma^{\gamma,s}_{c}(g,h) := \mu^{-\f12} Q^{\gamma,s}_{c}(\mu^{\f12}g,\mu^{\f12}h).
	\een
	The operator $\mathcal{L}^{\gamma,s}_{c}$ has a five dimensional kernel 
	\ben \label{null-based-on-N-limit-0}
	\ker \mathcal{L}^{\gamma,s}_{c} \colonequals     \mathrm{span} \{ \mu^{\f{1}{2}}, \mu^{\f{1}{2}}v_{1}, \mu^{\f{1}{2}}v_{2}, \mu^{\f{1}{2}}v_{3}, \mu^{\f{1}{2}}|v|^{2} \}.
	\een
	Let $\mathbb{P}_{0}$ be the projection onto the kernel. The self-adjoint operator $\mathcal{L}^{\gamma,s}_{c}$ has a coercive property in the orthogonal space of  $\ker \mathcal{L}^{\gamma,s}_{c}$ which is the counterpart of Boltzmann's H-Theorem
	in the linear level.
	More precisely,  it holds that 
	\ben \label{classical-coercivity-result-rough-1}
	\langle \mathcal{L}^{\gamma,s}_{c} f , f \rangle \sim |f-\mathbb{P}_{0}f|_{X}^{2}
	\een
	for some anisotropic norm $X$. Here $X = \mathcal{N}^{s,\gamma}$ in \cite{gressman2011global} and is the triple norm $|||\cdot|||$  in \cite{alexandre2011boltzmann} and \cite{alexandre2012boltzmann}. Both  norms 
	are defined through some integral so they  are not explicit.  An equivalent norm $|\cdot|_{\mathcal{L}^{s}_{\gamma/2}}$ will be explicitly used in the following. The nonlinear term $\Gamma^{\gamma,s}_{c}$ can be controlled by the norm $X$ as 
	$$|\langle \Gamma^{\gamma,s}_{c}(g, h), f \rangle| \lesssim |g|_{L^{2}}  |h|_{X} |f|_{X}.$$
	 By using this upper bound estimate and the coercivity estimate
	  \eqref{classical-coercivity-result-rough-1},  the global in time well-posedness in the perturbative framework can be established.

%
%
%
%
%
%
%
%
%
%
%
%

	For later use, we  now recall an explicit norm $|\cdot|_{\mathcal{L}^{s}_{\gamma/2}}$ which is equivalent to $\mathcal{N}^{s,\gamma}$ in \cite{gressman2011global} and $|||\cdot|||$ in \cite{alexandre2011boltzmann} and \cite{alexandre2012boltzmann}. 
	Firstly, 
	for $n \in \mathbb{N}, l \in \mathbb{R}$ and a function $f(v)$ on $\mathbb{R}^{3}$, denote
	\ben \label{Sobolev-norm}
	|f|_{H^{n}_{l}}^{2} \colonequals    \sum_{|\beta| \leq n} |W_{l} \partial_{\beta}f|_{L^{2}}^{2}, \quad |f|_{L^{2}_{l}}\colonequals   |f|_{H^{0}_{l}}, \quad |f|_{L^{\infty}_{l}} \colonequals   |W_{l}f|_{L^{\infty}}.
	\een
 Here $W_{l}(v) \colonequals (1+|v|^{2})^{l/2}$ for  $l \in \mathbb{R}$ is a weight function.
 In general, the Sobolev norm can be defined through Fourier transform. Denote by $\hat{f}$ the Fourier transform of $f$.
For $n, l \in \mathbb{R}$, set
	\ben \label{Sobolev-regularity-norm-real-index}
	|f|_{H^{n}}^{2} \colonequals    \int_{\mathbb{R}^{3}} (1+|\xi|^{2})^{n} |\hat{f}(\xi)|^{2} \mathrm{d}\xi, \quad |f|_{H^{n}_{l}}^{2}\colonequals   |W_{l}f|_{H^{n}}^{2}.
	\een

	Now we introduce an order-$n$ anisotropic norm $A^{n}$ through the projection on the
	real spherical harmonics denoted by  $\{ Y_{l}^{m} \}_{l \geq 0, -l\leq m \leq l}$, where 
	$(-\varDelta_{\mathbb{S}^2})Y_{l}^{m}=l(l+1)Y_{l}^{m}$ with $-\varDelta_{\mathbb{S}^2}$ being the 	Laplace-Beltrami operator on the unit sphere $\mathbb{S}^2$. 
For $n \in \mathbb{R}$, define the norm $A^{n}$ by
	\ben \label{anisotropic-regularity-norm-real-index}
	|f|_{A^{n}}^{2}\colonequals    \sum_{l=0}^\infty\sum_{m=-l}^{l} \int_{0}^{\infty}  \big( 1+l(l+1) \big)^{n} (f^{m}_{l}(r))^{2} r^{2} \mathrm{d}r, \een
	where $ f^{m}_{l}(r) = \int Y^{m}_{l}(\sigma) f(r \sigma) \mathrm{d}\sigma$ is the weighted integral of $f$ over the ball $B_{r}  := \{v \in \mathbb{R}^{3} : |v| =r\}$ with weight $Y^{m}_{l}$. The notation $A$ refers to ``anisotropic regularity''.
  Note that for a radial function $g(v)=g(|v|)$,  it holds 
	\ben \label{radial-L-infty-out}
	|g f|_{A^{n}} \leq |g|_{L^{\infty}} |f|_{A^{n}}.
	\een
	
Then  for $l \in \mathbb{R}, 0<s<1$,  the norm $\mathcal{L}^{s}_{l}$ is defined by 
	\ben \label{definition-of-norm-L-epsilon-gamma}
	|f|_{\mathcal{L}^{s}_{l}}^{2} \colonequals    |W_{l}f|^{2}_{L^{2}_{s}} + |W_{l}f|_{H^{s}}^{2} + |W_{l}f|_{A^{s}}^{2}.
	\een
	Recalling \eqref{Sobolev-norm}, \eqref{Sobolev-regularity-norm-real-index} and \eqref{anisotropic-regularity-norm-real-index}, the three norms on the right-hand of \eqref{definition-of-norm-L-epsilon-gamma} share a common weight function $W_{s}(r) = (1+|r|^{2})^{\frac{s}{2}}$ for $|v|, |\xi|, (l(l+1))^{\f12}$. 
	This norm as a summation of three terms  is general and essential for characterization of the linearized Boltzmann operator, see He-Zhou \cite{he2021boltzmann,he2022asymptotic} and Duan-He-Yang-Zhou \cite{duan2022solutions} for discussion in different settings.

	With the norm defined in \eqref{definition-of-norm-L-epsilon-gamma}, we can rewrite 
	\eqref{classical-coercivity-result-rough-1} as the following Lemma.
	\begin{lem}\label{classical-linearized-operator-dissipation} Let $-3< \gamma \leq 1, 0<s<1$. There are two generic  constants $\lambda_{*}, C_{*}$(depending on $\gamma$ and $s$) such that
		\ben \label{classical-coercivity-result}
		\lambda_{*} |(\mathbb{I}-\mathbb{P}_{0})f|_{\mathcal{L}^{s}_{\gamma/2}}^{2} \leq  \langle \mathcal{L}^{\gamma,s}_{c} f , f \rangle \leq C_{*}|(\mathbb{I}-\mathbb{P}_{0})f|_{\mathcal{L}^{s}_{\gamma/2}}^{2}.
		\een
	\end{lem}
	The proof of
	Lemma \ref{classical-linearized-operator-dissipation} consists of two steps. The first one is a rough  estimate as
		\ben \label{classical-coercivity-rough}
	  \langle \mathcal{L}^{\gamma,s}_{c} f , f \rangle + |f|_{L^{2}_{\gamma/2}}^{2}
	\sim |f|_{\mathcal{L}^{s}_{\gamma/2}}^{2}.
	\een
	Noting that without projection there is a lower order term $|f|_{L^{2}_{\gamma/2}}^{2}$  on the left hand side. 
	This  estimate can be  obtained by using the triple norm  $|||\cdot|||$ introduced in \cite{alexandre2012boltzmann}, cf.  Section 2 of \cite{alexandre2012boltzmann}.
	See also \cite{gressman2011global} using the norm $\mathcal{N}^{s,\gamma}$ and a more recent work \cite{alexandre2019global} using the norm $\mathcal{L}^{s}_{\gamma/2}$
	by different methods. The second step is the  spectral gap estimate 
			\ben \label{classical-coercivity-spectral-gap}
	\langle \mathcal{L}^{\gamma,s}_{c} f , f \rangle \gtrsim |(\mathbb{I}-\mathbb{P}_{0})f|_{L^{2}_{\gamma/2}}^{2},
	\een
   which is used to deal with the lower order term $|f|_{L^{2}_{\gamma/2}}^{2}$ to obtain the lower bound in \eqref{classical-coercivity-result}. The estimate \eqref{classical-coercivity-spectral-gap} is first proved in \cite{wang1952propagation} for the Maxwellian molecules $\gamma=0$. For general case, it was investigated by Pao in  \cite{pao1974boltzmann1,pao1974boltzmann2} using some technical tools from pseudo-differential
   operator theory. Different tools and
   constructive proofs were given in the elegant 
   papers \cite{baranger2005explicit,mouhot2006explicit,mouhot2007spectral}.
   See \cite{lerner2013phase} from the perspective of functional calculus.
   Recently, some new developments were made in \cite{he2021boltzmann,duan2022solutions} for the study of grazing limit to the Landau equation.

%
%
%
%
%

	\subsection{Equilibrium and the linearized equation} Temperature plays an important role in the study of quantum Boltzmann equation. For example, for Bose-Einsein particles, BEC happens in low temperature. We now recall  some basic facts  about temperature in the quantum context.
	Let us consider a homogeneous density $f=f(v) \geq 0$ with zero mean $\int v f(v) \mathrm{d}v=0$. For $k \geq 0$, the $k$-th moment is defined by
	\beno
	M_{k}(f) \colonequals   \int |v|^{k} f(v)  \mathrm{d}v.
	\eeno
	Denote  $M_{0} = M_{0}(f), M_{2}=M_{2}(f)$ for simplicity, and  $m$ be the mass of a particle. Then
	$m M_{0}$ and $\f{1}{2}m M_{2}$ are the total mass and kinetic energy per unit space volume.
	By referring to \cite{lu2004isotropic},
	the kinetic temperature $\bar{T}$ and the critical temperature $\bar{T}_{c}$ 
	 are defined by
	\ben \label{kinetic-temperature-and-critical}
	\bar{T} = \frac{1}{3 k_{B}} \frac{m M_{2}}{M_{0}}, \quad \bar{T}_{c} = \frac{m\zeta(5/2)}{2 \pi k_{B} \zeta(3/2)}
	\big( \frac{M_{0} }{\zeta(3/2)} \big)^{\f23},
	\een
	where $k_{B}$ is the Boltzmann constant and $\zeta(s) = \sum_{n=1}^{\infty} \frac{1}{n^{s}}$ is the Riemann zeta function.
	
	The equilibrium of the classical Boltzmann equation is the Maxwellian distribution with density
	function  
	\beno
	\mu_{\rho, v_{0}, T}(v) \colonequals   \rho T^{-\frac{3}{2}}e^{-\frac{1}{2T}|v-v_{0}|^{2}}, \quad \rho, T>0, v_{0} \in \mathbb{R}^{3}.
	\eeno
	Note that we drop the usual constant $(2\pi)^{-\frac{3}{2}}$ for simplicity. Now $(2\pi)^{\frac{3}{2}}\rho$ is density, $v_{0}$ is mean velocity  and $T$ is temperature. The  Bose-Einstein distribution has density 
	function
	\ben \label{equilibrium}
	\mathcal{M}_{\rho, v_{0}, T} \colonequals   \frac{\mu_{\rho, v_{0}, T}}{1 - \mu_{\rho, v_{0}, T}}, \quad \rho T^{-\frac{3}{2}} \leq 1. \een
	Now $\rho, v_{0}$ and $T$ do not represent density, mean velocity and temperature anymore, but are only three parameters. The ratio $\bar{T}/\bar{T}_{c}$ quantifies high and low temperature.
	In high temperature  $\bar{T}/\bar{T}_{c} > 1$,  the equilibrium of BBE equation is the Bose-Einstein distribution \eqref{equilibrium} with
	$\rho T^{-\frac{3}{2}} < 1$. In low temperature  $\bar{T}/\bar{T}_{c} < 1$,  the equilibrium of BBE equation is the Bose-Einstein distribution \eqref{equilibrium} with
	$\rho T^{-\frac{3}{2}} = 1$ plus a Dirac delta function. 
	 In the critical case  when $\bar{T}/\bar{T}_{c} = 1$, the equilibrium is \eqref{equilibrium} with
	$\rho T^{-\frac{3}{2}} = 1$. One can refer to \cite{lu2005boltzmann} for details.
	
	Without  loss of generality, we set $v_{0}=0, T=1$. 
	Then by using  $\rho$ to denote the fugacity, the equilibrium $\mathcal{M}_{\rho, v_{0}, T}$  can be written as 
	\ben \label{equilibrium-rho}
	\mathcal{M}_{\rho} \colonequals   \frac{\rho \mu}{1 - \rho \mu}, \quad 0 \leq \rho \leq 1. \een

   We are interested the following Cauchy problem on the inhomogeneous BBE equation
	\ben \label{quantum-Boltzmann-CP}
	\partial _t F +  v \cdot \nabla_{x} F=Q(F,F), ~~t > 0, x \in \Omega, v \in \R^3 ; \quad
	F|_{t=0}(x,v) = F_{0}(x,v),
	\een
	where the operator $Q$ is defined by \eqref{U-U-operator-scaling} with $B_{\phi}$ replaced by $B=B^{\gamma,s}$ given in \eqref{inverse-power-law-kernel}. For general initial datum $F_{0}$, it is 
	difficult to prove global well-posedness of the above Cauchy problem. 
	However,  we can consider the case when the  initial datum $F_{0}$ is close to an equilibrium $\mathcal{M}_{\rho}$. Recall that the solution of \eqref{quantum-Boltzmann-CP} conserves mass, momentum and energy. That is, for any $t \geq 0$,
	 \ben \label{conversation-mass-momentum-energy}
	 \int (1, v, |v|^{2}) F(t,x,v) \mathrm{d}x \mathrm{d}v = \int (1, v, |v|^{2}) F_{0}(x,v) \mathrm{d}x \mathrm{d}v.
	 \een
	 From now on, we assume 
	 \ben \label{F0-determines-equilibrium}
	 \int (1, v, |v|^{2}) \left(F_{0}(x,v) -  \mathcal{M}_{\rho}(v)\right) \mathrm{d}x \mathrm{d}v = 0.
	 \een
	
	To study the Cauchy problem, set 
	\ben \label{expansion-m-n}
	F = \mathcal{M}_{\rho} + \mathcal{N}_{\rho} f, \een
	where
	\ben \label{multiplier-rho}
	\mathcal{N}_{\rho} \colonequals   \frac{\rho^{\f{1}{2}} \mu^{\f{1}{2}}}{1 - \rho \mu}.  \een
	This special choice of $\mathcal{N}_{\rho}$ leads to the fact  that the  linearized quantum Boltzmann operator $\mathcal{L}^{\rho}$ defined in  \eqref{linearized-quantum-Boltzmann-operator-UU} is self-adjoint (see \eqref{self-joint-operator-quantum-case}).

	For simplicity of notations, denote $\mathcal{M}\colonequals  \mathcal{M}_{\rho}$, $ \mathcal{N}\colonequals  \mathcal{N}_{\rho}$. With  $F = \mathcal{M} + \mathcal{N} f$,
	 \eqref{quantum-Boltzmann-CP} becomes
	\ben\label{linearized-quantum-Boltzmann-eq} \left\{ \begin{aligned}
		&\partial _t f +  v \cdot \nabla_{x} f + \mathcal{L}^{\rho}f = \Gamma_{2}^{\rho}(f,f) + \Gamma_{3}^{\rho}(f,f,f), ~~t > 0, x \in \Omega, v \in \R^3 ,\\
		&f|_{t=0} = f_{0} = \frac{(1-  \rho\mu)F_{0} - \rho\mu}{\rho^{\f{1}{2}} \mu^{\f{1}{2}}}.
	\end{aligned} \right.
	\een
	Here the linearized quantum Boltzmann operator $\mathcal{L}^{\rho}$ is defined by
	\ben \label{linearized-quantum-Boltzmann-operator-UU}
	(\mathcal{L}^{\rho}f)(v)  \colonequals   \int  B \mathcal{N}_{*} \mathcal{N}^{\prime} \mathcal{N}^{\prime}_{*} \mathrm{S}(\mathcal{N}^{-1}f)
	\mathrm{d}\sigma \mathrm{d}v_{*},
	\een
		where 
	\ben  \label{symmetry-operator}
	\mathrm{S}(g) \colonequals     g + g_{*} - g^{\prime} - g^{\prime}_{*}.
	\een
	The bilinear term $\Gamma_{2}^{\rho}(\cdot,\cdot)$ and the trilinear term $\Gamma_{3}^{\rho}(\cdot,\cdot,\cdot)$ are given by
	\ben \label{definition-Gamma-2-epsilon}
	\Gamma_{2}^{\rho}(g,h) &\colonequals&   \mathcal{N}^{-1}\int
	B \Pi_{2}(g,h)
	\mathrm{d}\sigma \mathrm{d}v_{*},
	\\ \label{definition-Gamma-3-epsilon}
	\Gamma_{3}^{\rho}(g,h,\varrho)  &\colonequals&     \mathcal{N}^{-1}\int
	B  \big( (\mathcal{N}g)_{*}^{\prime} (\mathcal{N}h)^{\prime} ((\mathcal{N}\varrho)_{*} + \mathcal{N}\varrho) - (\mathcal{N}g)_{*} (\mathcal{N}h) ((\mathcal{N}\varrho)_{*}^{\prime} + (\mathcal{N}\varrho)^{\prime}) \big) \mathrm{d}\sigma \mathrm{d}v_{*}. \quad
	\quad
	\een
	Here, 
	\ben   \label{definition-A-2}
	\Pi_{2}(g,h) &\colonequals  & (\mathcal{N}g)^{\prime}_{*}(\mathcal{N}h)^{\prime}
	- (\mathcal{N}g)_{*}(\mathcal{N}h)
	\\ \label{line-1} &&+(\mathcal{N}g)^{\prime}_{*}(\mathcal{N}h)^{\prime}(\mathcal{M} +  \mathcal{M}_{*})
	-(\mathcal{N}g)_{*}(\mathcal{N}h)(\mathcal{M}^{\prime} +  \mathcal{M}_{*}^{\prime})
	\\ \label{line-2} &&+  (\mathcal{N}g)_{*}(\mathcal{N}h)^{\prime}( \mathcal{M}^{\prime}_{*} -  \mathcal{M})
	- (\mathcal{N}g)_{*}^{\prime}(\mathcal{N}h)( \mathcal{M}_{*} -  \mathcal{M}^{\prime})	
	\\ \label{line-3} &&+ (\mathcal{N}g)^{\prime}(\mathcal{N}h)(\mathcal{M}^{\prime}_{*} - \mathcal{M}_{*}) + (\mathcal{N}g)^{\prime}_{*}(\mathcal{N}h)_{*}(\mathcal{M}^{\prime} - \mathcal{M}).
	\een
	Note that the three operators $\mathcal{L}^{\rho}$, $\Gamma_{2}^{\rho}(\cdot,\cdot)$ and $\Gamma_{3}^{\rho}(\cdot,\cdot,\cdot)$ depend on $\rho$ through $\mathcal{M} =  \mathcal{M}_{\rho}$,  $\mathcal{N} =  \mathcal{N}_{\rho}$.

	By symmetry, the self-adjoint property of $\mathcal{L}^{\rho}$ satisfies
	\ben \label{self-joint-operator-quantum-case}
	\langle \mathcal{L}^{\rho}g, h \rangle = \frac{1}{4} \int B \mathcal{N} \mathcal{N}_{*} \mathcal{N}^{\prime} \mathcal{N}^{\prime}_{*} \mathrm{S}( \mathcal{N}^{-1} g ) \mathrm{S}(\mathcal{N}^{-1} h) \mathrm{d}V = \langle g,  \mathcal{L}^{\rho}h \rangle.
	\een
	For simplicity, here and in the rest of the paper,
	we use the notation $\mathrm{d}V \colonequals \mathrm{d}\sigma \mathrm{d}v_{*} \mathrm{d}v$.
	Therefore,  $\langle \mathcal{L}^{\rho}f, f \rangle \geq 0$. Recall that the kernel space $\ker \mathcal{L}^{\rho}$ of $\mathcal{L}^{\rho}$ is 
	\ben \label{def-kernel-quantum}
	\ker \mathcal{L}^{\rho} \colonequals     \mathrm{span} \{ \mathcal{N}_{\rho}, \mathcal{N}_{\rho}v_{1}, \mathcal{N}_{\rho}v_{2}, \mathcal{N}_{\rho}v_{3}, \mathcal{N}_{\rho}|v|^{2} \}.
	\een
	And the following four statements are equivalent 
	\ben \label{inner-zero}
	\langle \mathcal{L}^{\rho}f, f \rangle = 0 \Leftrightarrow
	\mathrm{S}( \mathcal{N}^{-1}_{\rho} f ) = 0 \Leftrightarrow
 \mathcal{N}^{-1}_{\rho} f \text{ is a collision invariant } \Leftrightarrow
	 f \in 	\ker \mathcal{L}^{\rho}.
	\een
	For $\rho>0$,
	let $\mathbb{P}_{\rho}$ be the projection on the kernel space $\ker \mathcal{L}^{\rho}$. For its explicit definition, see \eqref{definition-of-e-pm-rho}. The subscript $\rho$ is used to indicate dependence on the fugacity $\rho$. When $\rho=0$, $\mathbb{P}_{\rho}= \mathbb{P}_{0}$ is just  the  projection operator onto the kernel  \eqref{null-based-on-N-limit-0}
	for the classical Boltzmann equation.
	
		\subsection{Main result} 
	In order to show the stability of \eqref{quantum-Boltzmann-CP} under a small initial perturbation $f_{0}$, a key ingredient is the  coercivity estimate of  the linearized quantum Boltzmann operator $\mathcal{L}^{\rho}$.  The coercivity estimate is more subtle when 
	$\rho \to 1^-$ due to the critical temperature for BEC.
	
	The main result of this paper on the  coercivity estimate of $\mathcal{L}^{\rho}$ for any $0 \leq \rho<1$ is stated as follows.
	
	\begin{thm}\label{main-theorem} Let $-3<\gamma \leq  0<s<1$.
		There are two generic  constants $\lambda_{0}, C_{0}>0$ depending only on $\gamma,s$ such that
		\ben \label{lower-and-upper-bound}
		\lambda_{0} (1-\rho)^{13 \times 2^{N_{\gamma,s}} -3} \rho |(\mathbb{I}-\mathbb{P}_{\rho})f|_{\mathcal{L}^{s}_{\gamma/2}}^{2} \leq 
		\langle \mathcal{L}^{\rho} f , f \rangle \leq C_{0}(1-\rho)^{-4}\rho|(\mathbb{I}-\mathbb{P}_{\rho})f|_{\mathcal{L}^{s}_{\gamma/2}}^{2},
		\een
		where $N_{\gamma,s} := -1 $ for $-2s \leq  \gamma \leq 0$; and $N_{\gamma,s} := [\f{-\gamma-2s}{s}]$ for $-3 < \gamma < -2s$.
	\end{thm}
	
	The constant $C_{0}>0$ is explicitly given in \eqref{defintion-of-constants-C0}. For the case $-2s \leq  \gamma \leq 0$,
	the constant $\lambda_{0}$ is explicitly given in \eqref{defintion-of-constants-lambda0}, and for 
	the case $-3 < \gamma < -2s$, it is given implicitly  the proof. 
	
	In the following, we will give a few remarks on Theorem \ref{main-theorem}.
	
	\begin{rmk} \label{the-same-norm}
    By Theorem \ref{main-theorem}, for any fixed $0<\rho<1$, we have $\langle \mathcal{L}^{\rho} f , f \rangle \sim |(\mathbb{I}-\mathbb{P}_{\rho})f|_{\mathcal{L}^{s}_{\gamma/2}}^{2}$.
	Hence, this together with Lemma \ref{classical-linearized-operator-dissipation} show that  both the quantum linearized Boltzmann operator and the classical one are characterized by the same norm $|\cdot|_{\mathcal{L}^{s}_{\gamma/2}}$ except for different kernel spaces. 
	\end{rmk}
		\begin{rmk} \label{coercivity-estimate-key}
	 The lower bound in \eqref{lower-and-upper-bound} is the  coercivity estimate of $\mathcal{L}^{\rho}$ which is the key in proving the stability of \eqref{quantum-Boltzmann-CP} around $\mathcal{M}_{\rho}$. The lower bound in \eqref{lower-and-upper-bound} validates the  stability for any temperature above the critical value. We will have a brief discussion on  this 
	 in the next subsection.
	\end{rmk}
  \begin{rmk} \label{rho-to-1} 
The estimate \eqref{lower-and-upper-bound} is not optimal in terms of the power of $1-\rho$. A sharp estimate like $\langle \mathcal{L}^{\rho} f , f \rangle \sim h(1-\rho)\rho|(\mathbb{I}-\mathbb{P}_{\rho})f|_{\mathcal{L}^{s}_{\gamma/2}}^{2}$ for some 
precise function $h(\cdot)$
	  	is desired but is challenging. 
	\end{rmk}

We now present the main ideas of the proof  of Theorem \ref{main-theorem} as follows.

	\begin{enumerate}
		\item In order to use the coercivity estimate of  the classical Boltzmann operator, we first
		give  a relation between the two Dirichlet forms $\langle \mathcal{L}^{\rho} f , f \rangle$ and $\langle \mathcal{L}_{c} f , f \rangle$ where $\mathcal{L}_{c} = \mathcal{L}^{\gamma,s}_{c}$.  Then for any function $f$ in the orthogonal space $(\ker \mathcal{L}^{\rho})^{\perp}$, we define
		 a corresponding function $\varPhi_f$ in the orthogonal space  
		$(\ker \mathcal{L}_{c})^{\perp}$ such that the two Dirichlet forms $\langle \mathcal{L}^{\rho} f , f \rangle$ and $\langle \mathcal{L}_{c} \varPhi_f, \varPhi_f \rangle$ are mutually controlled by each other. This reveals the relationship between the two kernel spaces defined in  \eqref{Def-L-and-Gamma} and \eqref{def-kernel-quantum}, cf. 
		\eqref{def-w-f-psi-f}, \eqref{functions-in-the-spaces} and \eqref{key-relation-btw-quantum-and-classical} for details.  
		\item The above relation and the estimate in \eqref{classical-coercivity-result} for the classical operator $\mathcal{L}_{c}$ give the upper bound in \eqref{lower-and-upper-bound} for the quantum operator $\mathcal{L}^{\rho}$, and the lower bound up to 
		a  $L^{2}$ norm.
		 For hard potential when $\gamma+2s \ge 0$, since  $|f|_{\mathcal{L}^{s}_{\gamma/2}} \geq |f|_{L^{2}_{\gamma/2+s}} \geq |f|_{L^{2}}$,  the $L^{2}$ norm can be recovered by the orthogonal decompostion. 
			\item For soft potential $\gamma+2s<0$, motivated by \cite{alexandre2012boltzmann,he2021boltzmann} for Boltzmann equation, we apply an induction argument by starting from the case when $\gamma=-2s$ and using the gain of moment of order $s$ to control the $L^2$ norm. 
		Hence, for any $-3 <\gamma < -2s$,  we need $N_{\gamma,s} + 1= [\f{-\gamma-2s}{s}] + 1$ steps of  induction.   
	\end{enumerate}

	\subsection{Global well-posedness and stability}
	In this subsection, we will discuss the global well-posedness and stability of the BBE  equation \eqref{linearized-quantum-Boltzmann-eq}  in the perturbative framework and show its 
	dependence  on the fugacity parameter $\rho$.  
	If $\rho \ll 1$, then
	$M_{2}(\mathcal{M}_{\rho}) \sim \rho, M_{0}(\mathcal{M}_{\rho}) \sim \rho$
	and  by  \eqref{kinetic-temperature-and-critical} it holds that
	\beno
	\frac{\bar{T}}{\bar{T}_{c}} \sim \rho^{-\frac{2}{3}} \gg 1.
	\eeno
	Under this high temperature assumption, the global well-posedness of the spatially homogeneous BBE
	equation  with (slightly general than) hard sphere collision was firstly proved in \cite{li2019global}. In the inhomogeneous case, also for hard sphere collision, the global well-posedness was proved in
	\cite{bae2021relativistic} and  \cite{ouyang2022quantum}.  These two papers regard the fugacity $\rho$ as a constant without considering the  temperature effect. 
	For the non-cutoff kernel \eqref{inverse-power-law-kernel}, under the high temperature condition $\rho \ll 1$, the global well-posedness of the inhomogeneous BBE  equation \eqref{linearized-quantum-Boltzmann-eq} near $\mathcal{M}_{\rho}$ is proved in \cite{zhou2022global}.

  For  any $\rho\in (0,1)$, 
  with the key coercivity estimate in Theorem \ref{main-theorem}, we briefly
  discuss the  proof of the  global well-posedness of the Cauchy problem \eqref{linearized-quantum-Boltzmann-eq} and
   stability of $\mathcal{M}_{\rho}$. For this, we need to estimate the nonlinear terms $\Gamma_{2}^{\rho}(\cdot,\cdot)$ defined in \eqref{definition-Gamma-2-epsilon} and $\Gamma_{3}^{\rho}(\cdot,\cdot,\cdot)$ defined in \eqref{definition-Gamma-3-epsilon}. Such estimates are derived in \cite{zhou2022global} for $\rho \ll 1$. 
  Here we need to be careful when 
  $\rho \to 1^-$.
 For the bilinear term $\Gamma_{2}^{\rho}(\cdot,\cdot)$, following the proof of Theorem 3.2 in \cite{zhou2022global}, we can prove 
   \ben \label{gamma-2-upper-bound}
 |\langle \Gamma_{2}^{\rho}(g,h), f \rangle|
 \lesssim (1-\rho)^{-5} \rho^{\f12} \min\{|g|_{H^{2}}  |h|_{\mathcal{L}^{s}_{\gamma/2}} , |g|_{L^{2}}(|h|_{\mathcal{L}^{s}_{\gamma/2}}+| h|_{H^{s+2}_{\gamma/2}}+| h|_{H^{\frac{3}{2}}_{\gamma/2+s}})\} |f|_{\mathcal{L}^{s}_{\gamma/2}}.
 \een
 Comparing \eqref{gamma-2-upper-bound} with Theorem 3.2 in \cite{zhou2022global}, there is a new factor $(1-\rho)^{-5}$.
For the trilinear term $\Gamma_{3}^{\rho}(\cdot,\cdot,\cdot)$, following the proof of Theorem 4.1 in \cite{zhou2022global}, 
we can show that 
  \ben \label{rho-highest-derivative}
  |\langle \Gamma_{3}^{\rho}(g,h,\varrho), f\rangle| &\lesssim& (1-\rho)^{-4} \rho|g|_{H^{2}}|\mu^{\frac{1}{256}}h|_{H^{2}}|\varrho|_{\mathcal{L}^{s}_{\gamma/2}}|f|_{\mathcal{L}^{s}_{\gamma/2}}
  \\ \nonumber &&+ (1-\rho)^{-4}\rho|g|_{H^{3}}(|h|_{\mathcal{L}^{s}_{\gamma/2}}+|h|_{H^{s+2}_{\gamma/2}}
  +|h|_{H^{\frac{3}{2}}_{\gamma/2+s}})|\mu^{\frac{1}{256}}\varrho|_{L^{2}}|f|_{\mathcal{L}^{s}_{\gamma/2}},
  \\ \label{h-highest-derivative}
  |\langle \Gamma_{3}^{\rho}(g,h,\varrho), f\rangle| &\lesssim& (1-\rho)^{-4} \rho|g|_{H^{2}}|h|_{\mathcal{L}^{s}_{\gamma/2}}|\mu^{\frac{1}{256}}\varrho|_{H^{3}}|f|_{\mathcal{L}^{s}_{\gamma/2}},
  \\ \label{g-highest-derivative}
  |\langle \Gamma_{3}^{\rho}(g,h,\varrho), f\rangle| &\lesssim& (1-\rho)^{-4} \rho|g|_{L^{2}}(|h|_{\mathcal{L}^{s}_{\gamma/2}}+|h|_{H^{s+2}_{\gamma/2}}
  +|h|_{H^{\f{1}{2}}_{\gamma/2+s}})|\mu^{\frac{1}{256}}\varrho|_{H^{3}}|f|_{\mathcal{L}^{s}_{\gamma/2}}.
  \een

 Based on the coercivity estimate \eqref{lower-and-upper-bound}, the upper bound estimates \eqref{gamma-2-upper-bound}, \eqref{rho-highest-derivative}, \eqref{h-highest-derivative} and \eqref{g-highest-derivative}, we can expect that the macro-micro decomposition and nonlinear energy method leads to the proof of the global well-posedness   in some  suitable  function space. For hard potential when $\gamma+2s \geq 0$, we can work with
 the Sobolev  space by  using the energy functional
 \ben \label{definition-energy-functional}
 \mathcal{E}_{N}(f) \colonequals   \sum_{|\alpha|+|\beta| \leq N} \|\partial^{\alpha}_{x}\partial^{\beta}_{v} f \|_{L^{2}}^{2}.
 \een
 Here $L^{2} = L^{2}(\Omega \times \R^{3})$. In view of \eqref{gamma-2-upper-bound}, \eqref{rho-highest-derivative}, \eqref{h-highest-derivative} and \eqref{g-highest-derivative}, a certain order of regularity is necessary. As in \cite{zhou2022global}, let $N \geq 9$, we can derive the following  global well-posedness result for \eqref{linearized-quantum-Boltzmann-eq}. 	
  For some constant $M_{1}(N), M_{2}(N)>0$ depending only on $N$,  there exist a generic  constant $\delta>0$ such that 
 if
 \ben \label{condition-on-initial-data}
 F_{0} = \mathcal{M}_{\rho} + \mathcal{N}_{\rho}f_{0} \geq 0, \quad \mathcal{E}_{N}(f_{0}) \leq \delta \rho^{2N+1} (1-\rho)^{M_{1}(N)},
 \een
 then the Cauchy problem \eqref{linearized-quantum-Boltzmann-eq} has a unique global solution $f \in L^{\infty}([0,\infty);\mathcal{E}_{N})$ satisfying 
 \ben \label{uniform-estimate-global} F(t) \geq \mathcal{M}_{\rho} + \mathcal{N}_{\rho}f(t)  \geq 0, \quad  \sup_{t \geq 0}\mathcal{E}_{N}(f(t))  \leq C \rho^{-2N} (1-\rho)^{-M_{2}(N)} \mathcal{E}_{N}(f_{0}).\een
Note that for $0<\rho<1$, the global well-posedness is guaranteed if the initial perturbation is within a certain power of $\rho$ and $1-\rho$.

 The powers of $\rho$ and $1-\rho$ in \eqref{condition-on-initial-data} and \eqref{uniform-estimate-global} are needed because of   the coercivity estimate \eqref{lower-and-upper-bound} 
 \ben \label{lower-and-upper-bound-hard}
 \lambda_{0} (1-\rho)^{\f72} \rho |(\mathbb{I}-\mathbb{P}_{\rho})f|_{\mathcal{L}^{s}_{\gamma/2}}^{2} \leq 
 \langle \mathcal{L}^{\rho} f , f \rangle \leq C_{0}(1-\rho)^{-4}\rho|(\mathbb{I}-\mathbb{P}_{\rho})f|_{\mathcal{L}^{s}_{\gamma/2}}^{2}.
 \een
One can follow the estimation given in  \cite{zhou2022global} for fixed $\rho$ by 
inputting suitable powers of $\rho$ and $(1-\rho)$ in the energy estimates when needed.
 

 For  soft potential when $\gamma+2s < 0$,   
 the weighted Sobolev space is needed and  the energy functional for $l \geq -(\gamma+2s)N$ is 
 \ben \label{definition-energy-functional-2}
\mathcal{E}_{N,l}(f) \colonequals   \sum_{|\alpha|+|\beta| \leq N} \|W_{l+(\gamma+2s)|\beta|}\partial^{\alpha}_{x}\partial^{\beta}_{v} f \|_{L^{2}}^{2}.
 \een
 Here the weight function $W_{l+(\gamma+2s)|\beta|}$ is designed to deal with the convection term $v \cdot \nabla_{x} f$, cf.  \cite{guo2003classical}. Based on the coercivity estimate \eqref{lower-and-upper-bound} and upper bound estimates of 
 $|\langle \Gamma_{2}^{\rho}(g,h), W_{2l}f \rangle|$ and $|\langle \Gamma_{3}^{\rho}(g,h,\varrho), W_{2l}f\rangle|$(see  Theorem 5.3 and 5.6 in \cite{zhou2022global}), the proof of the global well-posedness of the Cauchy problem \eqref{linearized-quantum-Boltzmann-eq} can be expected with conclusion
 like \eqref{condition-on-initial-data} and \eqref{uniform-estimate-global}.
 Note that 
 the constants $M_{1}(N), M_{2}(N)$ 
 depend on $\gamma, s$. 
 

 
In order to focus on the coercivity of the linearized operator in this paper, we will not go into
the detailed estimation for the global well-posedness.

	\section{Proof of Theorem \ref{main-theorem}} \label{linear}
	We will present the  proof of Theorem \ref{main-theorem} in this section. The proof is divided into two parts. In the first subsection,  we will prove an upper bound for general $\gamma,s$ 
	and a  lower bound in the hard potential regime with $\gamma+2s \geq 0$. In the second subsection, we will prove the lower bound for the soft potential case when $\gamma+2s < 0$.

%
%
%

	\subsection{Lower bound for hard potential and upper bound}
	Recall $\mu(v) = \exp(-|v|^{2}/2)$ and set
	 \ben \label{equilibrium-rho-not-vanish}
	 M_{\rho} \colonequals  \frac{\mu}{1 - \rho \mu}, \quad N_{\rho} \colonequals    \frac{\mu^{\f{1}{2}}}{1 - \rho \mu}.
	 \een
	 By \eqref{equilibrium-rho} and \eqref{multiplier-rho} for 
	 the equilibrium $\mathcal{M}_{\rho}$ and multiplier $\mathcal{N}_{\rho}$, we have
	 \ben \label{M--with-mathcal-M}
	 \mathcal{M}_{\rho} = \rho M_{\rho}, \quad \mathcal{N}_{\rho} = \rho^{\f{1}{2}} N_{\rho}.
	 \een
	 By \eqref{linearized-quantum-Boltzmann-operator-UU}, since $\mathcal{N}_{\rho} = \rho^{\f{1}{2}} N_{\rho}$, for $\rho>0$, the null space of $\mathcal{L}^{\rho}$ is
	 \ben \label{null-based-on-N-not-vanish}
	 \ker \mathcal{L}^{\rho}  =  \mathrm{span} \{ N_{\rho}, N_{\rho}v_{1}, N_{\rho}v_{2}, N_{\rho}v_{3}, N_{\rho}|v|^{2} \}.
	 \een
	  
	 By \eqref{self-joint-operator-quantum-case} and \eqref{M--with-mathcal-M}, 
	 we have
	\ben \label{inner-of-quantum-operator}
	\langle \mathcal{L}^{\rho}f, f \rangle = \frac{\rho}{4} \int B N_{\rho} (N_{\rho})_{*} (N_{\rho})^{\prime} (N_{\rho})^{\prime}_{*}  \mathrm{S}^{2}( N_{\rho}^{-1} f )  \mathrm{d}V.
	\een

	Recalling \eqref{Def-L-and-Gamma}, the operator $\mathcal{L}_{c} = \mathcal{L}^{\gamma,s}_{c}$ reads
	\ben \label{linearized-classical-Boltzmann-operator}
	(\mathcal{L}_{c}f)(v)  =  \mu^{-1/2} \int  B \mu_{*} \mu \mathrm{S}(\mu^{-1/2}f)	\mathrm{d}\sigma \mathrm{d}v_{*},
	\een
	where $\mathrm{S}(\cdot)$ is defined by \eqref{symmetry-operator}.
	For later use, we
	define a functional $\mathcal{H}_{c} (\cdot)$ by
	\ben \label{definition-of-functional-H}
	\quad \mathcal{H}_{c} (f) :=  \langle \mathcal{L}_{c}f, f \rangle = \f{1}{4} \int B \mu \mu_{*}
	\mathrm{S}^{2}( \mu^{-1/2} f ) \mathrm{d}V.
	\een
	Note that the latter identity follows from the  symmetry property of the integral.
	
	We will compare the quantities in \eqref{inner-of-quantum-operator} and \eqref{definition-of-functional-H}. For this, let us define 
	\ben \label{definition-of-functional-J}
	\mathcal{J}_{\rho,\gamma,s} (f) := \frac{1}{4}   \int B^{\gamma,s} \mu \mu_{*}  \mathrm{S}^{2}( N_{\rho}^{-1} f )  \mathrm{d}V.
	\een
	Here we  specify the dependence of $\mathcal{J}_{\rho,\gamma,s}$ on $\gamma,s$
	through $B^{\gamma,s}$ defined in \eqref{inverse-power-law-kernel}, and on $\rho$ through $N_{\rho}$ defined in \eqref{M--with-mathcal-M}. For simplicity,  we sometimes use $\mathcal{J}_{\rho} $ for $\mathcal{J}_{\rho,\gamma,s}$. 
	
	If $0 \leq \rho < 1, 0< \mu \leq 1$, then
	\ben \label{M-rho-mu-N}
	\mu \leq M_{\rho} \leq (1-\rho)^{-1} \mu, \quad  \mu^{\f{1}{2}} \leq N_{\rho} \leq (1-\rho)^{-1} \mu^{\f{1}{2}}.
	\een
	As a consequence, by noting $\mu\mu_{*}=\mu^{\prime}\mu^{\prime}_{*}$, we have
	\ben \label{K-2-mu}
	\mu \mu_{*} \leq  N_{\rho} (N_{\rho})_{*} (N_{\rho})^{\prime} (N_{\rho})^{\prime}_{*} \leq (1-\rho)^{-4} \mu \mu_{*}.
	\een
	Hence, from \eqref{inner-of-quantum-operator} and \eqref{definition-of-functional-J}, it holds
	\ben \label{reduce-to-mu-mu-star}
	\rho \mathcal{J}_{\rho} (f) \leq \langle \mathcal{L}^{\rho}f, f \rangle \leq (1-\rho)^{-4} \rho \mathcal{J}_{\rho} (f).
	\een
	Observe $\mathcal{J}_{\rho} (f) = \mathcal{H}_{c} (\mu^{1/2} N_{\rho}^{-1} f)$. Thus, 
	 we can  use the estimate \eqref{classical-coercivity-result} 
	 to study $\mathcal{J}_{\rho} (f)$  by using  the connection between $\ker \mathcal{L}^{\rho}$ and $\ker \mathcal{L}_{c}$ in the following proof.	
	
	\begin{proof}[Proof of Theorem \ref{main-theorem}: ]
		In this part of the proof, we will show the upper bound estimate given in
		Theorem \ref{main-theorem}  and the lower bound for $\gamma+2s \geq 0$. It suffices to prove Theorem \ref{main-theorem} for $f \in (\ker \mathcal{L}^{\rho})^{\perp}$ because  $\langle \mathcal{L}^{\rho}f, f \rangle = \langle \mathcal{L}^{\rho}(\mathbb{I}-\mathbb{P}_{\rho})f, (\mathbb{I}-\mathbb{P}_{\rho})f \rangle$.

		For a function $f$, we define
		\ben \label{def-w-f-psi-f}
		w_{f} \colonequals N_{\rho}\mu^{-\f12} \mathbb{P}_{0} (N_{\rho}^{-1}\mu^{\f12}f) 
		, \quad \varPhi_{f}  \colonequals 
		(f - w_{f}) N_{\rho}^{-1} \mu^{1/2}.
		\een
		Recalling \eqref{null-based-on-N-not-vanish} and \eqref{null-based-on-N-limit-0},
		it is straightforward to check for $f \in (\ker \mathcal{L}^{\rho})^{\perp}$ that 
		\ben \label{functions-in-the-spaces}
		w_{f} \in \ker \mathcal{L}^{\rho}
		,  \quad \varPhi_{f} =  N_{\rho}^{-1}\mu^{\f12} f - \mathbb{P}_{0} (N_{\rho}^{-1}\mu^{\f12}f) \in (\ker \mathcal{L}_{c})^{\perp}.
		\een
		By the above construction, for $f \in (\ker \mathcal{L}^{\rho})^{\perp}$,
		we have the following key relation
		\ben \label{key-relation-btw-quantum-and-classical}
		\mathcal{J}_{\rho} (f) = \mathcal{J}_{\rho} (f - w_{f}) = \mathcal{J}_{\rho} (N_{\rho} \mu^{-1/2} \varPhi_{f})
		= \mathcal{H}_{c} (\varPhi_{f})
		\een

		By \eqref{key-relation-btw-quantum-and-classical},
		\eqref{functions-in-the-spaces}, \eqref{definition-of-functional-H} and  the upper bound in \eqref{classical-coercivity-result}, and by noting  $N_{\rho}^{-1} \mu^{1/2}=1-\rho \mu$,
		we have
		\ben \label{upper-bound-by-norm}
		\mathcal{J}_{\rho} (f) \leq C_{*} |\varPhi_{f}|_{\mathcal{L}^{s}_{\gamma/2}}^{2} 
		&\leq& 2 C_{*} (|(1-\rho \mu)f|_{\mathcal{L}^{s}_{\gamma/2}}^{2} + |\mathbb{P}_{0} ((1-\rho \mu)f)|_{\mathcal{L}^{s}_{\gamma/2}}^{2})
	\\
	&\leq &2 C_{*} (C_{1}+C_{2})|f|_{\mathcal{L}^{s}_{\gamma/2}}^{2},
		\een
		where $C_{1}, C_{2}$ are some generic constants such that for any $f$, 
		\ben \label{C1-inequality}
		|(1-\rho \mu)f|_{\mathcal{L}^{s}_{\gamma/2}}^{2} \leq C_{1}|f|_{\mathcal{L}^{s}_{\gamma/2}}^{2},
		\\\label{C2-inequality}
		|\mathbb{P}_{0} ((1-\rho \mu)f)|_{\mathcal{L}^{s}_{\gamma/2}}^{2} \leq C_{2}|f|_{\mathcal{L}^{s}_{\gamma/2}}^{2}.
		\een
		In fact, since $|\mu|_{H^{2}} \lesssim 1$, it holds 
		\beno
		|(1-\rho \mu)f|_{\mathcal{L}^{s}_{\gamma/2}} \lesssim |f|_{\mathcal{L}^{s}_{\gamma/2}} +
		|\mu f|_{\mathcal{L}^{s}_{\gamma/2}} \lesssim |f|_{\mathcal{L}^{s}_{\gamma/2}},
		\eeno
		which gives \eqref{C1-inequality}.  As $\mathbb{P}_{0}$ is the projection operator, it holds that
		\beno
		|\mathbb{P}_{0} ((1-\rho \mu)f)|_{\mathcal{L}^{s}_{\gamma/2}} \lesssim |\mu^{\f14}(1-\rho \mu)f|_{L^{2}} \lesssim |\mu^{\f14}f|_{L^{2}} \lesssim |f|_{\mathcal{L}^{s}_{\gamma/2}},
		\eeno
		which gives \eqref{C2-inequality}.
		Combining  \eqref{upper-bound-by-norm} and \eqref{reduce-to-mu-mu-star}, we obtain the upper bound in \eqref{lower-and-upper-bound} with
		\ben \label{defintion-of-constants-C0}
		C_{0} = 2 C_{*} (C_{1}+C_{2}).
		\een

		By \eqref{key-relation-btw-quantum-and-classical},
		\eqref{functions-in-the-spaces}, \eqref{definition-of-functional-H} and  the lower bound in \eqref{classical-coercivity-result}, we have
		\ben \label{J-rho-lower-bound}
		\mathcal{J}_{\rho} (f) \geq \lambda_{*} |\varPhi_{f}|_{\mathcal{L}^{s}_{\gamma/2}}^{2} \geq C_{3} (1-\rho)^{2} \lambda_{*} |N_{\rho} \mu^{-1/2} \varPhi_{f}|_{\mathcal{L}^{s}_{\gamma/2}}^{2} =  C_{3}\lambda_{*} (1-\rho)^{2}|f - w_{f}|_{\mathcal{L}^{s}_{\gamma/2}}^{2}.
		\een
		Here $C_{3}$ is  a generic  constant such that  
		\ben \label{norm-factor-out}
		|f|_{\mathcal{L}^{s}_{\gamma/2}}^{2} \geq C_{3}(1-\rho)^{2}|N_{\rho} \mu^{-1/2}f|_{\mathcal{L}^{s}_{\gamma/2}}^{2}. 
		\een
		For clear presentation,  the proof of \eqref{norm-factor-out} is postponed to Lemma \ref{factor-out}.

		We now derive that $|w_{f}|_{\mathcal{L}^{s}_{\gamma/2}}$ is bounded by $|\mu^{\f14}f|_{L^{2}}$ so that it is a lower order term. 
		Note that  $w_{f} = (a_{f} + b_{f} \cdot v + c_{f} |v|^{2}) N_{\rho}$ where $a_{f}, b_{f}, c_{f}$ are the constants given by
		\beno
		a_{f} = \f{5}{2} \langle f, N_{\rho}^{-1} \mu \rangle - \f{1}{2} \langle f, |v|^{2}N_{\rho}^{-1} \mu \rangle, \quad
		b_{f} = \langle f, v N_{\rho}^{-1} \mu \rangle, \quad
		c_{f} =   \f{1}{6} \langle f, |v|^{2}N_{\rho}^{-1} \mu \rangle - \f{1}{2} \langle f, N_{\rho}^{-1} \mu \rangle.
		\eeno
		Obviously, 
		\ben \label{w-f-upper-bound}
		|w_{f}|_{\mathcal{L}^{s}_{\gamma/2}} \lesssim (|a_{f}| + |b_{f}| + |c_{f}|) |\mu^{-\f18}N_{\rho}|_{\mathcal{L}^{s}_{\gamma/2}}
		\lesssim (1-\rho)^{-\f34}|\mu^{\f{1}{4}}f|_{L^{2}},
		\een
		where we have used Lemma \ref{N-rho-norm-L2-H1} (to be proved later) to get 
		\ben \label{upper-bound-of-N-rho}
		|(1+|v|^{2})N_{\rho}|_{\mathcal{L}^{s}_{\gamma/2}} \lesssim |\mu^{-\f18}N_{\rho}|_{\mathcal{L}^{s}_{\gamma/2}} \lesssim
		|\mu^{-\f14}N_{\rho}|_{H^{1}} \sim (1-\rho)^{-\f34}.
		\een
		By \eqref{w-f-upper-bound}, let $C_{4}$ be the generic constant such that $|w_{f}|_{\mathcal{L}^{s}_{\gamma/2}}^{2} \leq C_{4} (1-\rho)^{-\f32}|\mu^{\f14}f|_{L^{2}}^{2}$. Then we get
		\ben \label{leading-with-l2-missing}
		|f - w_{f}|_{\mathcal{L}^{s}_{\gamma/2}}^{2} \geq \f{1}{2}|f|_{\mathcal{L}^{s}_{\gamma/2}}^{2} - |w_{f}|_{\mathcal{L}^{s}_{\gamma/2}}^{2}
		\geq \f{1}{2}|f|_{\mathcal{L}^{s}_{\gamma/2}}^{2} - C_{4} (1-\rho)^{-\f32} |\mu^{\f14}f|_{L^{2}}^{2}.
		\een

		We now assume $\gamma+2s \geq 0$.
		Then by \eqref{definition-of-norm-L-epsilon-gamma},  $|f|_{\mathcal{L}^{s}_{\gamma/2}} \geq |f|_{L^{2}_{\gamma/2+s}} \geq |f|_{L^{2}}$. For $f \in (\ker \mathcal{L}^{\rho})^{\perp}$,
		 \eqref{functions-in-the-spaces} implies that $f \perp w_{f}$.  Hence, 
		\ben \label{orthoganal}
		|f - w_{f}|_{\mathcal{L}^{s}_{\gamma/2}}^{2} \geq |f - w_{f}|_{L^{2}}^{2} = |f|_{L^{2}}^{2} + |w_{f}|_{L^{2}}^{2} \geq |f|_{L^{2}}^{2}.
		\een
		Suitable combination of \eqref{orthoganal} and \eqref{leading-with-l2-missing} gives
		\beno
		|f - w_{f}|_{\mathcal{L}^{s}_{\gamma/2}}^{2} \geq \f{1}{2(1+ C_{4}(1-\rho)^{-\f32})}|f|_{\mathcal{L}^{s}_{\gamma/2}}^{2}.
		\eeno
		The above estimate together with \eqref{J-rho-lower-bound} yield
		\ben \label{lower-bound-J-rho-gamma-s}
		\mathcal{J}_{\rho} (f) \geq   \f{C_{3}\lambda_{*}(1-\rho)^{2}}{2(1+ C_{4}(1-\rho)^{-\f32})} |f|_{\mathcal{L}^{s}_{\gamma/2}}^{2} \geq \f{C_{3}\lambda_{*}(1-\rho)^{\f72}}{2(1+ C_{4})} |f|_{\mathcal{L}^{s}_{\gamma/2}}^{2}.
		\een
		By \eqref{reduce-to-mu-mu-star}, 
		 we prove the lower bound in \eqref{lower-and-upper-bound} with
		\ben \label{defintion-of-constants-lambda0}
		\lambda_{0}=\f{C_{3}\lambda_{*}}{2(1+ C_{4})}.
		\een
		And this completes the first part of the proof.
	\end{proof}

In the following, we prove the two estimates used in the above proof in  two lemmas.
Firstly, 
we prove \eqref{norm-factor-out} in the following Lemma.
	\begin{lem}\label{factor-out}
	\ben  \label{singular-factor-out}
	|\f{1}{1-\rho \mu} f|_{\mathcal{L}^{s}_{\gamma/2}}  \lesssim (1-\rho)^{-1} |f|_{\mathcal{L}^{s}_{\gamma/2}}.
	\een
\end{lem}
\begin{proof} Let $g=\f{1}{1-\rho \mu}, h= W_{\gamma/2}f$.
	By \eqref{definition-of-norm-L-epsilon-gamma}, the norm $\mathcal{L}^{s}_{\gamma/2}$ has three parts
	\ben \label{product-norm-L-epsilon-gamma}
	|\f{1}{1-\rho \mu} f|_{\mathcal{L}^{s}_{\gamma/2}}^{2} =  |g h|^{2}_{L^{2}_{s}} + |g h|_{H^{s}}^{2} + |g h|_{A^{s}}^{2}.
	\een
	Note that
	\ben \label{L-infty-singularity}
	|g|_{L^{\infty}} = (1-\rho)^{-1}.
	\een
	For the first and third parts in \eqref{product-norm-L-epsilon-gamma}, by \eqref{L-infty-singularity} and \eqref{radial-L-infty-out}, we have
	\beno
	|g h|^{2}_{L^{2}_{s}} + |g h|_{A^{s}}^{2} \leq
	(1-\rho)^{-2}|W_{\gamma/2}f|^{2}_{L^{2}_{s}} + (1-\rho)^{-2}|W_{\gamma/2}f|_{A^{s}}^{2}.
	\eeno
	For the fractional Sobolev norm $H^{s}$, we claim
	\ben \label{fractional-Sobolev-out}
	|g h|_{H^{s}} \lesssim (1-\rho)^{-1} |h|_{H^{s}},
	\een
	so that we have \eqref{singular-factor-out}.
	
	Now it remains to prove \eqref{fractional-Sobolev-out}. Note that
	\beno
	|g h|_{H^{s}}^{2} \sim \int \f{|g(x)h(x)-g(y)h(y)|^{2}}{|x-y|^{3+2s}} \mathrm{d}x \mathrm{d}y + |g h|_{L^{2}}^{2}.
	\eeno
	Therefore, it suffices to derive
	\beno
	\mathcal{I} :=\int \f{|g(x)h(x)-g(y)h(y)|^{2}}{|x-y|^{3+2s}} \mathrm{d}x \mathrm{d}y \lesssim (1-\rho)^{-2} |h|_{H^{s}}^{2}.
	\eeno
	Note that in the proof  all the integrals are taken with respect to variables in  $\R^{3}$. 
	
	If $|x-y| \geq 1$, by using \eqref{L-infty-singularity}, we have
	\beno
	\int \mathrm{1}_{|x-y| \geq 1}\f{|g(x)h(x)-g(y)h(y)|^{2}}{|x-y|^{3+2s}} \mathrm{d}x \mathrm{d}y \lesssim |g h|_{L^{2}}^{2}
	\lesssim (1-\rho)^{-2}|h|_{L^{2}}^{2}.
	\eeno
	When $|x-y| \leq 1$, we write
	\beno
	&& \int \mathrm{1}_{|x-y| \leq 1}\f{|g(x)h(x)-g(y)h(y)|^{2}}{|x-y|^{3+2s}} \mathrm{d}x \mathrm{d}y 
	\\ &\lesssim&
	\int \mathrm{1}_{|x-y| \leq 1}\f{g^{2}(x)|h(x)-h(y)|^{2}}{|x-y|^{3+2s}} \mathrm{d}x \mathrm{d}y
	+ \int \mathrm{1}_{|x-y| \leq 1}\f{h^{2}(y)|g(x)-g(y)|^{2}}{|x-y|^{3+2s}} \mathrm{d}x \mathrm{d}y.
	\eeno
	The first term in the last inequality  is bounded by
	\beno
	\int \mathrm{1}_{|x-y| \leq 1}\f{g^{2}(x)|h(x)-h(y)|^{2}}{|x-y|^{3+2s}} \mathrm{d}x \mathrm{d}y \lesssim (1-\rho)^{-2}|h|_{H^{s}}^{2}.
	\eeno
	Now it remains to prove
	\beno
	\mathcal{J} := \int \mathrm{1}_{|x-y| \leq 1}\f{h^{2}(y)|g(x)-g(y)|^{2}}{|x-y|^{3+2s}} \mathrm{d}x \mathrm{d}y \lesssim (1-\rho)^{-2}|h|_{H^{s}}^{2}.
	\eeno
	Let us use the parameter $\lambda>0$ to replace $0<\rho<1$ by $e^{-\lambda^{2}} = \rho$, then
	\beno
	g(x) = \left(1 - \exp(-\f{|x|^{2}}{2} - \lambda^{2})\right)^{-1}.
	\eeno
	Direct computation gives
	\ben \label{g-diff-expression}
	g(x) - g(y) = \f{\exp(- \lambda^{2})\left(\exp(-\f{|x|^{2}}{2})-\exp(-\f{|y|^{2}}{2})\right)}{\left(1 - \exp(-\f{|x|^{2}}{2} - \lambda^{2})\right)\left(1 - \exp(-\f{|y|^{2}}{2} - \lambda^{2})\right)}.
	\een
	If  $\lambda \geq \f12$ (equivalent to $\rho \leq e^{-1/4}$), then the denominator in \eqref{g-diff-expression} is bounded so that 
	\beno
	|g(x) - g(y)| \lesssim |\exp(-\f{|x|^{2}}{2})-\exp(-\f{|y|^{2}}{2})| \lesssim |x-y|
	\eeno
	which gives
	\beno
	\mathcal{J} \lesssim \int \mathrm{1}_{|x-y| \leq 1}\f{h^{2}(y)}{|x-y|^{1+2s}} \mathrm{d}x \mathrm{d}y \lesssim |h|_{L^{2}}^{2}
	\lesssim (1-\rho)^{-2}|h|_{L^{2}}^{2}.
	\eeno
	
	Now it remains to estimate $\mathcal{J}$ when  $0<\lambda \leq \f12$.
	For this, we divide $\mathcal{J}$ into three parts,
	\beno
	\mathcal{J} \leq \int \mathrm{1}_{|x-y| \leq 1, |x| \geq 1/2} (\cdots)
	+ \int \mathrm{1}_{|x-y| \leq 1, |y|\geq 1/2} (\cdots)
	+ \int \mathrm{1}_{|x| \leq 1/2, |y| \leq 1/2} (\cdots) := \mathcal{J}_{1} + \mathcal{J}_{2} + \mathcal{J}_{3}.
	\eeno
	For  $\mathcal{J}_{1}$, we have $1 - \exp(-\f{|x|^{2}}{2} - \lambda^{2}) \geq 1 - \exp(-\f{1}{8})$ and $1 - \exp(-\f{|y|^{2}}{2} - \lambda^{2}) \geq 1 - \exp(-\lambda^{2}) \sim \lambda^{2}$.
	Thus,  
	\beno
	|g(x) - g(y)| \lesssim \lambda^{-2} |\exp(-\f{|x|^{2}}{2})-\exp(-\f{|y|^{2}}{2})| \lesssim \lambda^{-2}|x-y|
	\eeno
	which gives
	\beno
	\mathcal{J}_{1} \lesssim \int \mathrm{1}_{|x-y| \leq 1}\f{h^{2}(y)}{|x-y|^{1+2s}} \mathrm{d}x \mathrm{d}y \lesssim \lambda^{-4}|h|_{L^{2}}^{2}.
	\eeno
	Similarly, we can show that  $\mathcal{J}_{2}$ is also bounded by  $\lambda^{-4}|h|_{L^{2}}^{2}$.
	
	We now turn to  $\mathcal{J}_{3}$ for which $|x| \leq 1/2, |y| \leq 1/2$.  Then   we have
	\beno
	1 - \exp(-\f{|x|^{2}}{2} - \lambda^{2}) \sim \f{|x|^{2}}{2} + \lambda^{2} \sim |x|^{2} + \lambda^{2},
	\\
	|\exp(-\f{|x|^{2}}{2})-\exp(-\f{|y|^{2}}{2})| = |\exp(-\f{|y|^{2}}{2}) \left(\exp(\f{|y|^{2}}{2}-\f{|x|^{2}}{2})-1\right)| \sim \left| |y|^{2} - |x|^{2} \right|.
	\eeno
	Hence, 
	\ben \label{phi-x-y-upper-bound}
	|g(x) - g(y)|  \sim \f{\left| |y|^{2} - |x|^{2} \right|}{(|x|^{2} + \lambda^{2})(|y|^{2} + \lambda^{2})} := \phi(x,y) \lesssim \lambda^{-2}.
	\een
	For simplicity of notations, set
	\beno
	K(x,y) = \phi^{2}(x,y)h^{2}(y)|x-y|^{-3-2s}.
	\eeno
	We further divide $\mathcal{J}_{3}$ into three parts:
	\beno
	\mathcal{J}_{3} &\leq& \int \mathrm{1}_{|x| \leq 1/2, |y| \leq 1/2, |x-y| \leq \lambda} K(x,y) \mathrm{d}x \mathrm{d}y + \int \mathrm{1}_{|x| \leq 1/2, |y| \leq 1/2, |x-y| \geq \lambda, |x| \geq |y|} K(x,y) \mathrm{d}x \mathrm{d}y
	\\&&+ \int \mathrm{1}_{|x| \leq 1/2, |y| \leq 1/2, |x-y| \geq \lambda, |x| \leq |y|} K(x,y) \mathrm{d}x \mathrm{d}y
	:= \mathcal{K}_{1} + \mathcal{K}_{2} + \mathcal{K}_{3}.
	\eeno
	For $\mathcal{K}_{1}$, since  $|x-y| \leq \lambda$,  we use
	\beno
	\phi(x,y) \lesssim \f{|x-y| (|x|+|y|)}{(|x|^{2} + \lambda^{2})(|y|^{2} + \lambda^{2})}
	\eeno
	to get
	\beno
	\mathcal{K}_{1} \lesssim \int \mathrm{1}_{|x| \leq 1/2, |y| \leq 1/2, |x-y| \leq \lambda} \f{(|x|^{2}+|y|^{2}) h^{2}(y)|x-y|^{-1-2s}}{(|x|^{2} + \lambda^{2})^{2}(|y|^{2} + \lambda^{2})^{2}}
	\mathrm{d}x \mathrm{d}y
	\\ \lesssim \lambda^{-6+2s} \int \mathrm{1}_{|y| \leq 1/2, |x-y| \leq \lambda} \f{h^{2}(y)|x-y|^{-1-2s}}{|y|^{2s}}
	\mathrm{d}x \mathrm{d}y  \lesssim \lambda^{-4} |h|^{2}_{H^{s}}.
	\eeno
	Here we have used
	\beno
	\f{|x|^{2}+|y|^{2}}{(|x|^{2} + \lambda^{2})^{2}(|y|^{2} + \lambda^{2})^{2}} \lesssim \lambda^{-6+2s} |y|^{-2s},
	\eeno
	the estimate $\int \mathrm{1}_{|u| \leq \lambda} |u|^{-1-2s} \mathrm{d}u \lesssim \lambda^{2-2s}$ and the Hardy's inequality $\int \mathrm{1}_{|y| \leq 1/2} \f{h^{2}(y)}{|y|^{2s}} \mathrm{d}y \lesssim |h|^{2}_{H^{s}}$.
	
	For $\mathcal{K}_{2}$, since  $|y| \leq |x|$, we use
	\beno
	\phi(x,y) \lesssim \f{|x|^{2}}{(|x|^{2} + \lambda^{2})(|y|^{2} + \lambda^{2})} \lesssim \lambda^{-2+s} |y|^{-s}
	\eeno
	to get
	\beno
	\mathcal{K}_{2}  \lesssim \lambda^{-4+2s} \int \mathrm{1}_{|y| \leq 1/2, |x-y| \geq \lambda} \f{h^{2}(y)|x-y|^{-3-2s}}{|y|^{2s}}
	\mathrm{d}x \mathrm{d}y  \lesssim \lambda^{-4} |h|^{2}_{H^{s}},
	\eeno
	where we have used $\int \mathrm{1}_{|u| \geq \lambda} |u|^{-3-2s} \mathrm{d}u \lesssim \lambda^{-2s}$.
	
	For $\mathcal{K}_{3}$, we first divide it into two terms
	\beno
	\mathcal{K}_{3} \lesssim \int \mathrm{1}_{|x| \leq 1/2, |y| \leq 1/2, |x-y| \geq \lambda, |x| \leq |y|}
	\phi^{2}(x,y)h^{2}(x)|x-y|^{-3-2s} \mathrm{d}x \mathrm{d}y
	\\+ \int \mathrm{1}_{|x| \leq 1/2, |y| \leq 1/2, |x-y| \geq \lambda, |x| \leq |y|}
	\phi^{2}(x,y)|h(x)-h(y)|^{2}|x-y|^{-3-2s} \mathrm{d}x \mathrm{d}y.
	\eeno
	Note that the first term is exactly $\mathcal{K}_{2}$ that is bounded by $\lambda^{-4} |h|^{2}_{H^{s}}$.
	By \eqref{phi-x-y-upper-bound}, the second term is also bounded by $\lambda^{-4} |h|^{2}_{H^{s}}$. 
	In summary, we have
	\beno
	\mathcal{J} \lesssim \lambda^{-4} |h|^{2}_{H^{s}} \lesssim (1 - \rho)^{-2} |h|^{2}_{H^{s}},
	\eeno
	where we have used  $1 - \rho =1 - e^{-\lambda^{2}} \sim \lambda^{2}$ for $0<\lambda \leq \f12$. 
	And this completes the proof of the lemma.
\end{proof}

We now prove
the last inequality in \eqref{upper-bound-of-N-rho} in the following Lemma.
\begin{lem}\label{N-rho-norm-L2-H1} Let $\f{1}{4} \leq a \leq 1$, then it holds that
	\beno
	|\mu^{-\f12+a}N_{\rho}|_{L^{2}} \sim (1-\rho)^{-\f14}, \quad |\nabla_{v}(\mu^{-\f12+a}N_{\rho})|_{L^{2}} \sim (1-\rho)^{-\f34}.
	\eeno
\end{lem}
\begin{proof} By \eqref{equilibrium-rho-not-vanish},
	\ben \label{L2-expression}
	|\mu^{-\f12+a}N_{\rho}|_{L^{2}}^{2} = \int \f{\mu^{2a}}{(1-\rho \mu)^{2}} \mathrm{d}v.
	\een
	Again, denote $\rho=e^{-\lambda^{2}}$. Then
	\beno
	1-\rho \mu = 1 - \exp(-\f{|v|^{2}}{2} - \lambda^{2}).
	\eeno
	If $\lambda \geq \f12$ (equivalent to $\rho \leq e^{-1/4}$), then 
	$(1-\rho \mu)^{-1} \sim 1$ so that
	\beno
	|\mu^{-\f12+a}N_{\rho}|_{L^{2}}^{2} \sim 1 \sim (1-\rho)^{-\f12}.
	\eeno
	
	For  $\lambda < \f12$,  noting that $1- \rho = 1- e^{-\lambda^{2}} \sim \lambda^{2}$, we divide the integral into two parts,
	\ben \label{L2-expression-2-parts}
	|\mu^{-\f12+a}N_{\rho}|_{L^{2}}^{2} = \int \mathrm{1}_{|v| \geq 1} 
	\f{\mu^{2a}}{(1-\rho \mu)^{2}}  \mathrm{d}v + \int \mathrm{1}_{|v| \leq 1}\f{\mu^{2a}}{(1-\rho \mu)^{2}} \mathrm{d}v.
	\een
	For the first part, $(1-\rho \mu)^{-1} \sim 1$ so that
	\beno
	\int \mathrm{1}_{|v| \geq 1}
	\f{\mu^{2a}}{(1-\rho \mu)^{2}}  \mathrm{d}v \sim 1.
	\eeno
	For the second part, since $\f{|v|^{2}}{2} + \lambda^{2} \leq 1$, then $1-\rho \mu = 1-\exp(-\f{|v|^{2}}{2} - \lambda^{2}) \sim |v|^{2}+\lambda^{2}$ and 
	\beno
	\int \mathrm{1}_{|v| \leq 1} \f{\mu^{2a}}{(1-\rho \mu)^{2}} \mathrm{d}v \sim 
	\int \mathrm{1}_{|v| \leq 1}\f{1}{(|v|^{2}+\lambda^{2})^{2}} \mathrm{d}v =
	\lambda^{-1}\int \mathrm{1}_{|v| \leq \lambda^{-1}}\f{1}{(|v|^{2}+1)^{2}} \mathrm{d}v \sim \lambda^{-1}.
	\eeno
	Hence, we have 
	\beno
	|\mu^{-\f12+a}N_{\rho}|_{L^{2}}^{2} \sim \lambda^{-1} \sim (1-\rho)^{-\f12}.
	\eeno
	
	Since $\nabla_{v}\mu = - v \mu$, direct computation gives 
	\beno
	\nabla_{v} (\mu^{-\f12+a}N_{\rho}) = (1-\rho \mu)^{-1} \nabla_{v} \mu^{a} +  \mu^{a} \nabla_{v} (1-\rho \mu)^{-1} 
	\\ = (1-\rho \mu)^{-1} (- a \mu^{a} v)  -  \mu^{a} (1-\rho \mu)^{-2}    v \rho \mu = - (1-\rho \mu)^{-2}(a + (1-a) \rho \mu) \mu^{a} v.
	\eeno
	Then 
	\beno
	|\nabla_{v} (\mu^{-\f12+a}N_{\rho})|_{L^{2}}^{2} = \int (1-\rho \mu)^{-4}(a + (1-a) \rho \mu)^{2} \mu^{2a} |v|^{2}  \mathrm{d}v 
	\sim \int (1-\rho \mu)^{-4}  \mu^{2a} |v|^{2}  \mathrm{d}v.
	\eeno
	
	Similarly, we can show that 
	\beno
	|\nabla_{v} (\mu^{-\f12+a}N_{\rho})|_{L^{2}}^{2} \sim (1-\rho)^{-\f32}.
	\eeno
	In fact, if  $\lambda \geq \f12$(equivalent to $\rho \leq e^{-1/4}$), then $(1-\rho \mu)^{-1} \sim 1$ and 
	\beno
	|\nabla_{v} (\mu^{-\f12+a}N_{\rho})|_{L^{2}}^{2} \sim  \int  \mu^{2a} |v|^{2}  \mathrm{d}v \sim 1 \sim (1-\rho)^{-\f12}.
	\eeno
	
	For $\lambda < \f12$, note that  $1- \rho = 1- e^{-\lambda^{2}} \sim \lambda^{2}$. We divide the integral into two parts,
	\ben \label{H1-expression-2-parts}
	|\nabla_{v} (\mu^{-\f12+a}N_{\rho})|_{L^{2}}^{2} \sim \int \mathrm{1}_{|v| \geq 1} (1-\rho \mu)^{-4} \mu^{2a} |v|^{2} \mathrm{d}v + \int \mathrm{1}_{|v| \leq 1} (1-\rho \mu)^{-4} \mu^{2a} |v|^{2} \mathrm{d}v.
	\een
	For the first part, since $(1-\rho \mu)^{-1} \sim 1$,  we have
	\beno
	\int \mathrm{1}_{|v| \geq 1} (1-\rho \mu)^{-4} \mu^{2a} |v|^{2} \mathrm{d}v \sim 1.
	\eeno
	For the second part, since $\f{|v|^{2}}{2} + \lambda^{2} \leq 1$, then $1-\exp(-\f{|v|^{2}}{2} - \lambda^{2}) \sim |v|^{2}+\lambda^{2}$ and 
	\beno
	\int \mathrm{1}_{|v| \leq 1} (1-\rho \mu)^{-4} \mu^{2a} |v|^{2} \mathrm{d}v &\sim&
	\int \mathrm{1}_{|v| \leq 1}\f{1}{(|v|^{2}+\lambda^{2})^{4}} |v|^{2} \mathrm{d}v 
	\\&=&
	\lambda^{-3} \int \mathrm{1}_{|v| \leq \lambda^{-1}}\f{1}{(|v|^{2}+1)^{2}} |v|^{2} \mathrm{d}v \sim \lambda^{-3}.
	\eeno
	Hence, 
	\beno
	|\nabla_{v} (\mu^{-\f12+a}N_{\rho})|_{L^{2}}^{2} \sim \lambda^{-3} \sim (1-\rho)^{-\f32}.
	\eeno
	And this completes the proof of the lemma.
\end{proof}

	\subsection{Lower bound for soft potential}
For soft potential with $\gamma+2s < 0$, the first inequality in \eqref{orthoganal} does not hold. 
Motivated by  \cite{alexandre2012boltzmann,he2021boltzmann} about exchanging the 
kinetic component  in the cross-section with a weight of velocity, we will apply an induction argument based on the estimate for the 
hard potential obtained in the previous subsection and the gain of moment of order $s$. As the first step, we reduce the case when $-3s\le \gamma<-2s$  to  $\gamma=-2s$, and then by induction  to the whole range $-3<\gamma <-2s$. For this, we first introduce a  weight function 
	\ben\label{specialweightfun} U_{\delta}(v) := W(\delta v) = (1+ \delta^{2}|v|^{2})^{1/2} \geq \max\{\delta|v|,1\}. \een
	Here $\delta$ is a sufficiently small parameter to be chosen later. We now give two lemmas on some integrals involving $U_{\delta}$.
	
	\begin{lem}\label{difference-term-complication} Let $\alpha, \beta<0<s, \delta <1$ with $\alpha+2s>-3$.
		Set 
		\beno 
		X(\beta,\delta) := \delta^{-\beta}\left((U^{\beta/2}_{\delta})^{\prime}(U^{\beta/2}_{\delta})^{\prime}_{*}-
		U^{\beta/2}_{\delta}(U^{\beta/2}_{\delta})_{*}\right)^{2}, \quad 
		b_{s}(\theta)  :=
		\sin^{-2-2s}\frac{\theta}{2} \mathrm{1}_{0 \leq \theta \leq \pi/2}.
		\eeno
		Then for $v \in \R^{3}$, 
		\ben \label{chi-W-difference-part} \int b_{s}(\theta) |v-v_{*}|^{\alpha}
		X(\beta,\delta) \mu_{*} \mathrm{d}\sigma \mathrm{d}v_{*}
		\lesssim  \max\{\frac{1}{s}, \frac{1}{1-s}\} \frac{|\beta|^{2s}\delta^{2s}}{\alpha+2s+3}  \langle v \rangle^{\alpha+\beta+2s}.\een
	\end{lem}
	\begin{proof}
		Note that 
		\beno X(\beta,\delta) \lesssim \delta^{-\beta}(U^{\beta}_{\delta})^{\prime}_{*}\left((U^{\beta/2}_{\delta})^{\prime}-(U^{\beta/2}_{\delta})\right)^{2}
		+\delta^{-\beta}U^{\beta}_{\delta}\left((U^{\beta/2}_{\delta})^{\prime}_{*}-(U^{\beta/2}_{\delta})_{*}\right)^{2}
		:= A_{1} + A_{2}. \eeno
		We only estimate $A_{1}$ because $A_{2}$ can be estimated similarly.
		
		For $a \leq 0$, one has
		\ben \label{derivative-bounds}
		|\nabla U^{a}_{\delta}| \lesssim  |a| \delta U^{a}_{\delta}
		\een
		which gives
		\beno \big((U^{a}_{\delta})^{\prime}-(U^{a}_{\delta})\big)^{2}  = \big|\int_{0}^{1} (\nabla U^{a}_{\delta})(v(\kappa))\cdot (v^{\prime}-v)\mathrm{d}\kappa\big|^{2}
		\lesssim a^{2} \delta^{2} \sin^{2}\f{\theta}{2}|v-v_{*}|^{2}  \int_{0}^{1}  U^{2a}_{\delta}(v(\kappa))\mathrm{d}\kappa,
		\eeno
		where $v(\kappa) = (1-\kappa)v + \kappa v^{\prime}$.
		Thanks to $|v^{\prime}_{*}|^{2}+|v(\kappa)|^{2} \sim |v|^{2}+|v_{*}|^{2}$, we have
		\ben \label{natural-bound-using-conservation}
		\delta^{-2a}(U^{2a}_{\delta})^{\prime}_{*} U^{2a}_{\delta}(v(\kappa)) \lesssim \langle v \rangle^{2a}
		\een 
		which implies
		\ben \label{order-2-out}
		\delta^{-2a}(U^{2a}_{\delta})^{\prime}_{*}\big((U^{a}_{\delta})^{\prime}-(U^{a}_{\delta})\big)^{2} \lesssim a^{2} \delta^{2} \sin^{2}\f{\theta}{2}|v-v_{*}|^{2} \langle v \rangle^{2a}.
		\een

		Divide the integral $ \int b_{s}(\theta) |v-v_{*}|^{\alpha}
		A_{1} \mu_{*} \mathrm{d}\sigma \mathrm{d}v_{*}$ into two parts: $\mathcal{I}_{\leq}$ and $\mathcal{I}_{\geq}$ corresponding to $|\beta|\delta|v-v_{*}|\leq 1$ and  $|\beta|\delta|v-v_{*}|\geq 1$.
		When $|\beta|\delta|v-v_{*}|\leq 1$, using \eqref{order-2-out} for $a = \beta/2$,
		we have
		\beno
		\mathcal{I}_{\leq} \lesssim   \beta^{2} \delta^{2} \langle v \rangle^{\beta} \int \mathrm{1}_{|v-v_{*}|\leq |\beta|^{-1}\delta^{-1}} b_{s}(\theta) \sin^{2}\f{\theta}{2} |v-v_{*}|^{2+\alpha} \mu_{*} \mathrm{d}\sigma \mathrm{d}v_{*}
		\lesssim \frac{1}{1-s}\frac{|\beta|^{2s}\delta^{2s}}{\alpha+2s+3} \langle v \rangle^{\alpha+\beta+2s}.
		\eeno
		
		When $|\beta|\delta|v-v_{*}|\geq 1$. We further divide  the integral $\mathcal{I}_{\geq}$ into two parts: 
		$\mathcal{I}_{\geq,\leq}$ and $\mathcal{I}_{\geq,\geq}$ corresponding to
		$\sin\f{\theta}{2} \leq |\beta|^{-1}\delta^{-1}|v-v_{*}|^{-1}$ and
		$\sin\f{\theta}{2} \geq |\beta|^{-1}\delta^{-1}|v-v_{*}|^{-1}$. Using \eqref{order-2-out} for $a = \beta/2$,
		we have
		\beno
		\mathcal{I}_{\geq, \leq} &\lesssim&  \beta^{2}  \delta^{2} \langle v \rangle^{\beta} \int \mathrm{1}_{|v-v_{*}| \geq |\beta|^{-1} \delta^{-1}} \mathrm{1}_{\sin\f{\theta}{2} \leq |\beta|^{-1} \delta^{-1}|v-v_{*}|^{-1}} b_{s}(\theta) \sin^{2}\f{\theta}{2} |v-v_{*}|^{2+\alpha} \mu_{*} \mathrm{d}\sigma \mathrm{d}v_{*}
		\\ &\lesssim& \frac{1}{1-s}\frac{|\beta|^{2s}\delta^{2s}}{\alpha+2s+3}  \langle v \rangle^{\alpha+\beta+2s}.
		\eeno
		
		For the remaining part with $\sin\f{\theta}{2} \geq |\beta|^{-1}\delta^{-1}|v-v_{*}|^{-1}$, it holds from \eqref{natural-bound-using-conservation} that
		 $A_{1} \lesssim \langle v \rangle^{\beta}$ and
		\beno
		\mathcal{I}_{\geq, \geq} &\lesssim&  \langle v \rangle^{\beta} \int \mathrm{1}_{|v-v_{*}|\geq |\beta|^{-1} \delta^{-1}} \mathrm{1}_{\sin\f{\theta}{2} \geq |\beta|^{-1}\delta^{-1}|v-v_{*}|^{-1}} b_{s}(\theta) |v-v_{*}|^{\alpha} \mu_{*} \mathrm{d}\sigma \mathrm{d}v_{*}
		\\ &\lesssim&  \frac{1}{s}\frac{|\beta|^{2s}\delta^{2s}}{\alpha+2s+3} \langle v \rangle^{\alpha+\beta+2s}.
		\eeno
		Combining the above estimates completes the proof of the lemma. 
	\end{proof}
	
	\begin{lem}\label{difference-term-complication-2} Let $\alpha, \beta<0<s, \delta <1$ with $\alpha+2s>-3$.
		Set $\varphi_{\beta,\delta}:= (1-U^{\beta/2}_{\delta})\mu^{\f{1}{2}}$, then
		\beno \int b_{s}(\theta) |v-v_{*}|^{\alpha}
		(\varphi_{\beta,\delta}^{\prime} -  \varphi_{\beta,\delta})^{2} \mathrm{d}\sigma \mathrm{d}v
		\lesssim \max\{\frac{\beta^{2}\delta^{2}}{s}, \frac{1}{1-s}\}  \frac{|\beta|^{2s}\delta^{2s}}{\alpha+2s+3}  \langle v_{*} \rangle^{\alpha+2s}.\eeno
	\end{lem}
	\begin{proof} By \eqref{derivative-bounds}, we get
		\ben \label{derivative-bounds-2}
		|\varphi_{\beta,\delta}| \lesssim  |\beta| \delta \mu^{\f{1}{4}}, \quad |\nabla \varphi_{\beta,\delta}| \lesssim  |\beta| \delta \mu^{\f{1}{4}}.
		\een
		
		Again divide the integral into two parts: $\mathcal{I}_{\leq}$ and $\mathcal{I}_{\geq}$ corresponding to $|\beta|\delta|v-v_{*}|\leq 1$ and  $|\beta|\delta|v-v_{*}|\geq 1$.
		When $|\beta|\delta|v-v_{*}|\leq 1$, by using the second estimate in \eqref{derivative-bounds-2} and the change of variable $v \to v(\kappa)$,
		we have
		\beno
		\mathcal{I}_{\leq} &\lesssim&   \beta^{2} \delta^{2}   \int \mathrm{1}_{|v-v_{*}|\leq |\beta|^{-1}\delta^{-1}} b_{s}(\theta) \sin^{2}\f{\theta}{2} |v-v_{*}|^{2+\alpha} \mu^{\f12}(v(\kappa)) \mathrm{d}\sigma \mathrm{d}v \mathrm{d}\kappa
		\\&\lesssim&  \frac{1}{1-s}\frac{|\beta|^{2s}\delta^{2s}}{\alpha+2s+3} \langle v_{*} \rangle^{\alpha+2s}.
		\eeno
		
		When $|\beta|\delta|v-v_{*}|\geq 1$, we also  divide the  integral $\mathcal{I}_{\geq}$ into two parts: 
		$\mathcal{I}_{\geq,\leq}$ and $\mathcal{I}_{\geq,\geq}$ corresponding to
		$\sin\f{\theta}{2} \leq |\beta|^{-1}\delta^{-1}|v-v_{*}|^{-1}$ and
		$\sin\f{\theta}{2} \geq |\beta|^{-1}\delta^{-1}|v-v_{*}|^{-1}$.
		For $\mathcal{I}_{\geq,\leq}$, by using the second estimate in \eqref{derivative-bounds-2} and the change of variable $v \to v(\kappa)$,
		we have
		\beno
		\mathcal{I}_{\geq, \leq} &\lesssim&   \beta^{2}\delta^{2}   \int \mathrm{1}_{|v-v_{*}|\geq |\beta|^{-1} \delta^{-1}} \mathrm{1}_{\sin\f{\theta}{2} \leq |\beta|^{-1}\delta^{-1}|v-v_{*}|^{-1}} b_{s}(\theta) \sin^{2}\f{\theta}{2} |v-v_{*}|^{2+\alpha} \mu^{\f12}(v(\kappa)) \mathrm{d}\sigma \mathrm{d}v \mathrm{d}\kappa
		\\ &\lesssim&  \frac{1}{1-s}\frac{|\beta|^{2s}\delta^{2s}}{\alpha+2s+3} \langle v_{*} \rangle^{\alpha+2s}.
		\eeno
		
		For the remaining part with $\sin\f{\theta}{2} \geq |\beta|^{-1}\delta^{-1}|v-v_{*}|^{-1}$,  by using the first estimate in \eqref{derivative-bounds-2} and the change of variable $v \to v^{\prime}$,
		we have
		\beno
		\mathcal{I}_{\geq, \geq} &\lesssim&   \beta^{2}\delta^{2} \int \mathrm{1}_{|v-v_{*}|\geq |\beta|^{-1}\delta^{-1}} \mathrm{1}_{\sin\f{\theta}{2} \geq |\beta|^{-1}\delta^{-1}|v-v_{*}|^{-1}} b_{s}(\theta) |v-v_{*}|^{\alpha} (\mu^{\f12} + (\mu^{\f12})^{\prime}) \mathrm{d}\sigma \mathrm{d}v
		\\ &\lesssim&  \frac{1}{s}\frac{|\beta|^{2s+2}\delta^{2s+2}}{\alpha+2s+3} \langle v_{*} \rangle^{\alpha+2s}.
		\eeno
		Combining the above estimates completes the proof of the lemma.
	\end{proof}
	
	We now give explicitly the definition of the projection operator $\mathbb{P}_{\rho}$.
	Recalling \eqref{null-based-on-N-not-vanish},
	we construct an orthogonal basis for $\ker \mathcal{L}^{\rho} $ 
	 as follows
	\ben \label{definition-of-di}
	\{ d^{\rho}_{i} \}_{1 \leq i \leq 5} \colonequals    \{ N_{\rho}, N_{\rho} v_{1}, N_{\rho} v_{2}, N_{\rho} v_{3}, N_{\rho}|v|^{2} -
	\langle N_{\rho}|v|^{2} , N_{\rho}\rangle|N_{\rho}|^{-2}_{L^{2}}N_{\rho} \}.
	\een
	Note that $\langle N_{\rho}|v|^{2} , N_{\rho}\rangle|N_{\rho}|^{-2}_{L^{2}}N_{\rho}$ is the projection of $N_{\rho}|v|^{2}$ on $N_{\rho}$.
	By normalizing $\{ d^{\rho}_{i} \}_{1 \leq i \leq 5}$, an orthonormal basis of $\ker \mathcal{L}^{\rho} $ can be obtained as
	\ben \label{definition-of-e-pm-rho}
	\{ e^{\rho}_{i} \}_{1 \leq i \leq 5} \colonequals    \{ \frac{d^{\rho}_{i}}{|d^{\rho}_{i}|_{L^{2}}} \}_{1 \leq i \leq 5}.
	\een
	With this orthonormal basis, the projection operator $\mathbb{P}_{\rho}$ on the null space  $\ker \mathcal{L}^{\rho} $ is defined by
	\ben \label{definition-projection-operator}
	\mathbb{P}_{\rho}f \colonequals    \sum_{i=1}^{5} \langle f, e^{\rho}_{i}\rangle e^{\rho}_{i}.
	\een

	We are now ready to  prove the coercivity estimate of $\mathcal{L}^{\rho}$ when $-3 < \gamma < -2s$ by induction.

	\begin{proof}[Proof of Theorem \ref{main-theorem}, continued:] In this part, we will show
		the  lower bound in \eqref{lower-and-upper-bound} in the case of $\gamma < -2s$. Again, it suffices to consider $f \in (\ker \mathcal{L}^{\rho})^{\perp}$ since  $\langle \mathcal{L}^{\rho}f, f \rangle = \langle \mathcal{L}^{\rho}(\mathbb{I}-\mathbb{P}_{\rho})f, (\mathbb{I}-\mathbb{P}_{\rho})f \rangle$.
		Note that for the case of $\gamma+2s \geq 0$, we already have the coercivity estimate.
		
		For simplicity of notations, we set
		\beno \mathbb{A}(f, g):=(f_{*}g + f g_{*} - f^{\prime}_{*}g^{\prime} - f^{\prime} g^{\prime}_{*}), \quad\mathbb{F}(f, g):= \mathbb{A}^{2}(f, g).\eeno
		 Then
		$\mathcal{J}_{\rho,\gamma,s}(f) = \f{1}{4}\int B^{\gamma,s} \mathbb{F}(\mu^{\f12}, g) \mathrm{d}V$ where $g = (1-\rho \mu)f$.
		
	We divide the proof  into four steps.
		
		\noindent{\it Step 1: Localization of $\mathcal{J}_{\rho,\gamma,s}(f)$.}
		By  \eqref{specialweightfun} and if $a \le 0$,  we get \beno |v - v_{*}|^{-a}
		\leq  C_{a} \delta^{a} ((\delta|v|)^{-a}+(\delta|v_{*}|)^{-a})
		\leq 2C_{a} \delta^{a} U^{-a}_{\delta}(v)U^{-a}_{\delta}(v_{*})  \eeno
		which gives
		\beno |v - v_{*}|^{a} \gtrsim \delta^{-a} U^{a}_{\delta}(v)U^{a}_{\delta}(v_{*}). \eeno
		With $a=\gamma+2s$, we have
		\beno \mathcal{J}_{\rho,\gamma,s}(f) \gtrsim \delta^{-\gamma-2s}\int b_{s}(\theta)|v-v_{*}|^{-2s}  U^{\gamma+2s}_{\delta}(U^{\gamma+2s}_{\delta})_{*}\mathbb{F}(\mu^{\f12}, g) \mathrm{d}V.
		\eeno
		By setting  $h= U^{\gamma/2+s}_{\delta}, \phi=\mu^{\f12}$ and commuting the weight function $U^{\gamma+2s}_{\delta}(U^{\gamma+2s}_{\delta})_{*}$ with $\mathbb{F}(\cdot,\cdot)$, we have
		\ben \label{move-inside} \nonumber
		&& U^{\gamma+2s}_{\delta}(U^{\gamma+2s}_{\delta})_{*}\mathbb{F}(\mu^{\f12}, g)  
		\\  &=&
		h^{2}_{*}h^{2}\mathbb{F}(\phi, g)
		=
		\left(h h_{*}\left(\phi_{*}g + \phi g_{*}\right) -
		h h_{*}\left(\phi^{\prime}_{*}g^{\prime} + \phi^{\prime} g^{\prime}_{*}\right)\right)^{2}
		\nonumber \\&=& \left(h h_{*}\left(\phi_{*}g + \phi g_{*}\right) -
		h^{\prime} h^{\prime}_{*}\left(\phi^{\prime}_{*}g^{\prime} + \phi^{\prime} g^{\prime}_{*}\right)
		+ \left(h^{\prime} h^{\prime}_{*}-
		h h_{*}\right)
		\left(\phi^{\prime}_{*}g^{\prime} + \phi^{\prime} g^{\prime}_{*}\right)
		\right)^{2}
		\nonumber \\  &\geq&  \frac{1}{2}  \left(h h_{*}\left(\phi_{*}g + \phi g_{*}\right) -
		h^{\prime} h^{\prime}_{*}\left(\phi^{\prime}_{*}g^{\prime} + \phi^{\prime} g^{\prime}_{*}\right) \right)^{2}\nonumber -\left(h^{\prime} h^{\prime}_{*}-
		h h_{*}\right)^{2}
		\left(\phi^{\prime}_{*}g^{\prime} + \phi^{\prime} g^{\prime}_{*}\right)^{2}\nonumber
		\\  &=&  \frac{1}{2}  \mathbb{F}(h \phi, h g)
		-\left(h^{\prime} h^{\prime}_{*}-
		h h_{*}\right)^{2}
		\left(\phi^{\prime}_{*}g^{\prime} + \phi^{\prime} g^{\prime}_{*}\right)^{2}.
		\een
		Thus, 
		\ben \label{move-inside-sym} \mathcal{J}_{\rho,\gamma,s}(f) &\gtrsim&
		\frac{1}{2}\delta^{-\gamma-2s}
		\int b_{s}(\theta)|v-v_{*}|^{-2s}  \mathbb{F}(U^{\gamma/2+s}_{\delta}\mu^{\f12}, U^{\gamma/2+s}_{\delta}g) \mathrm{d}V
		\\  && - \delta^{-\gamma-2s}\int b_{s}(\theta)|v-v_{*}|^{-2s} \left(h^{\prime} h^{\prime}_{*}-
		h h_{*}\right)^{2}
		\left(\phi^{\prime}_{*}g^{\prime} + \phi^{\prime} g^{\prime}_{*}\right)^{2} \mathrm{d}V. \nonumber
		\een
		We further rewrite 
		 $\mathbb{F}(U^{\gamma/2+s}_{\delta}\mu^{\f12}, U^{\gamma/2+s}_{\delta}g)$ 
		as  $\mathbb{F}(\mu^{\f12}, U^{\gamma/2+s}_{\delta}g)$ plus some correction terms. That is,
		\ben \label{move-outside-asym} \mathbb{F}(U^{\gamma/2+s}_{\delta}\mu^{\f12}, U^{\gamma/2+s}_{\delta}g)  &=& \mathbb{A}^{2}(U^{\gamma/2+s}_{\delta}\mu^{\f12}, U^{\gamma/2+s}_{\delta}g)
		\nonumber \\&=&\left( \mathbb{A}(\mu^{\f12}, U^{\gamma/2+s}_{\delta}g)  - \mathbb{A}\big((1-U^{\gamma/2+s}_{\delta})\mu^{\f12}, U^{\gamma/2+s}_{\delta}g\big) \right)^{2}
		\nonumber \\&\geq&\frac{1}{2} \mathbb{A}^{2}(\mu^{\f12}, U^{\gamma/2+s}_{\delta}g)
		-\mathbb{A}^{2}\big((1-U^{\gamma/2+s}_{\delta})\mu^{\f12}, U^{\gamma/2+s}_{\delta}g\big)
		\nonumber \\&=&\frac{1}{2}\mathbb{F}(\mu^{\f12}, U^{\gamma/2+s}_{\delta}g) -  \mathbb{F}((1-U^{\gamma/2+s}_{\delta})\mu^{\f12}, U^{\gamma/2+s}_{\delta}g).
		\een
		By symmetry and noting $\phi=\mu^{\f12}$, 
		we have
		\ben \nonumber && \int b_{s}(\theta)|v-v_{*}|^{-2s} \left(h^{\prime} h^{\prime}_{*}-
		h h_{*}\right)^{2}
		\left(\phi^{\prime}_{*}g^{\prime} + \phi^{\prime} g^{\prime}_{*}\right)^{2} \mathrm{d}V 
		\\ \label{sym-term} &\leq& 4\int b_{s}(\theta)|v-v_{*}|^{-2s} \left(h^{\prime} h^{\prime}_{*}-
		h h_{*}\right)^{2}\mu_{*}g^{2} \mathrm{d}V.
		\een
		By \eqref{move-inside-sym}, \eqref{move-outside-asym} and \eqref{sym-term}, we get
		\ben \label{separate-key-part} \mathcal{J}_{\rho,\gamma,s}(f) &\gtrsim&
		\frac{1}{4}\delta^{-\gamma-2s}
		\int b_{s}(\theta)|v-v_{*}|^{-2s} \mathbb{F}(\mu^{\f12}, U^{\gamma/2+s}_{\delta}g) \mathrm{d}V
		\\  && - \frac{1}{2}\delta^{-\gamma-2s}\int b_{s}(\theta)|v-v_{*}|^{-2s} \mathbb{F}((1-U^{\gamma/2+s}_{\delta})\mu^{\f12}, U^{\gamma/2+s}_{\delta}g) \mathrm{d}V
		\nonumber
		\\  && - 4\delta^{-\gamma-2s}\int b_{s}(\theta)|v-v_{*}|^{-2s} \left(h^{\prime} h^{\prime}_{*}-
		h h_{*}\right)^{2} \mu_{*}g^{2} \mathrm{d}V
		:= \frac{1}{4}J_{1} - \frac{1}{2}J_{2} -  4J_{3}. \nonumber
		\een
		
		\noindent{\it Step 2: Estimates of $J_i(i=1,2,3)$.} We will give the estimates term by term.
		
		{\it \underline{Lower bound of $J_1$.}}
		We claim that when $\delta$ is suitably small,
		\ben \label{estimate-J-1}  J_{1} \gtrsim \delta^{-\gamma-2s} (1-\rho)^{\f72} |f|^{2}_{L^{2}_{\gamma/2+s}}.  \een
		By \eqref{lower-bound-J-rho-gamma-s} and \eqref{defintion-of-constants-lambda0}, we have
		\ben \label{eta-pure-lower-bound} \mathcal{J}_{\rho,-2s,s}(F) \geq \lambda_{0} (1-\rho)^{\f72}|(\mathbb{I}-\mathbb{P}_{\rho})F|^{2}_{\mathcal{L}^{s}_{-s}}
		\geq \lambda_{0} (1-\rho)^{\f72} |(\mathbb{I}-\mathbb{P}_{\rho})F|^{2}_{L^{2}}. \een
		Applying \eqref{eta-pure-lower-bound} with $F=U^{\gamma/2+s}_{\delta}f$
		and using $(a-b)^{2} \geq a^{2}/2-b^{2}$, we have
		\beno  J_{1}&=&\delta^{-\gamma-2s}
		\int b_{s}(\theta)|v-v_{*}|^{-2s} \mathbb{F}(\mu^{\f12}, U^{\gamma/2+s}_{\delta}g) \mathrm{d}V
		=4 \delta^{-\gamma-2s} \mathcal{J}_{\rho,-2s,s}(U^{\gamma/2+s}_{\delta}f)
		\\&\gtrsim& \delta^{-\gamma-2s} (1-\rho)^{\f72} |(\mathbb{I}-\mathbb{P}_{\rho})(U^{\gamma/2+s}_{\delta}f)|^{2}_{L^{2}}.
		\eeno
		To get \eqref{estimate-J-1},   it suffices to show that for $a \leq 0$, the estimate 
		\ben \label{delta-small-projection-to-not-projection}
		|(\mathbb{I}-\mathbb{P}_{\rho})(U^{a}_{\delta}f)|^{2}_{L^{2}} \gtrsim |f|^{2}_{L^{2}_{a}}
		\een
		holds if $\delta$ is small enough. First, using $(a-b)^{2} \geq a^{2}/2-b^{2}$, we get
		\beno 
		|(\mathbb{I}-\mathbb{P}_{\rho})(U^{a}_{\delta}f)|^{2}_{L^{2}} \geq  \frac{1}{2} |U^{a}_{\delta}f|^{2}_{L^{2}} - |\mathbb{P}_{\rho}(U^{a}_{\delta}f)|^{2}_{L^{2}}
		\eeno
		Since $\delta \leq 1$ and $a \leq 0$,  $U^{a}_{\delta} \geq W_{a}$. Hence, 
		\ben \label{leading-term} |U^{a}_{\delta}f|^{2}_{L^{2}} \geq |f|^{2}_{L^{2}_{a}}.\een
		
		We now estimate $|\mathbb{P}_{\rho}(U^{a}_{\delta}f)|_{L^{2}}$ for $f \in (\ker \mathcal{L}^{\rho})^{\perp}$. Recalling \eqref{definition-projection-operator} for the definition of $\mathbb{P}_{\rho}$
		and  by the condition $\mathbb{P}_{\rho}f = 0$, we have
		\beno  \mathbb{P}_{\rho}(U^{a}_{\delta}f) = \sum_{i=1}^{5} e^{\rho}_{i}\int e^{\rho}_{i}U^{a}_{\delta}f \mathrm{d}v
		= \sum_{i=1}^{5} e^{\rho}_{i}\int e^{\rho}_{i}(U^{a}_{\delta}-1)f \mathrm{d}v .\eeno
		By \eqref{derivative-bounds}, for $a \leq 0$, 
		\ben\label{lclfact1} 0 \leq  1-U^{a}_{\delta}(v) = U^{a}_{\delta}(0)-U^{a}_{\delta}(v)\lesssim |a|\delta |v|.\een
		Thus
		$$ \big|\int e^{\rho}_{i}(U^{a}_{\delta}-1)f \mathrm{d}v\big| \lesssim |a| \delta |\mu^{\f18}f|_{L^{2}},$$
		where  $|v|$ cancels the singularity of $e^{\rho}_{i}$ at $v=0$.
		Therefore, 
		\ben \label{large-velocity-part-12}
		|\mathbb{P}_{\rho}(U^{a}_{\delta}f)|^{2}_{L^{2}} 
		 \lesssim a^{2} \delta^{2} |\mu^{\f18}f|^{2}_{L^{2}}.\een
		By combining the estimates \eqref{leading-term} and \eqref{large-velocity-part-12}
		and choosing $\delta$ suitably small, 
		we obtain \eqref{delta-small-projection-to-not-projection}.

		{\it \underline{Upper bound of $J_2$.}}
		For simplicity of notations, set $\varphi_{\gamma,\delta}= (1-U^{\gamma/2+s}_{\delta})\mu^{\f12}, \psi_{\gamma,\delta}=U^{\gamma/2+s}_{\delta}(1-\rho \mu)f$. Then
		\ben  \nonumber
		J_{2} &=&
		\delta^{-\gamma-2s}\int b_{s}(\theta)|v-v_{*}|^{-2s} \mathbb{F}((1-U^{\gamma/2+s}_{\delta})\mu^{\f12}, U^{\gamma/2+s}_{\delta}(1-\rho \mu)f) \mathrm{d}V
		\\ \nonumber &=&\delta^{-\gamma-2s}\int b_{s}(\theta)|v-v_{*}|^{-2s} \mathbb{F}(\varphi_{\gamma,\delta}, \psi_{\gamma,\delta}) \mathrm{d}V
		\\ \label{J2-part} &\lesssim& \delta^{-\gamma-2s}\int b_{s}(\theta)|v-v_{*}|^{-2s} (\varphi_{\gamma,\delta}^{2})_{*}(\psi_{\gamma,\delta}^{\prime}-\psi_{\gamma,\delta})^{2} \mathrm{d}V
		\\ \nonumber &&+ \delta^{-\gamma-2s}\int b_{s}(\theta)|v-v_{*}|^{-2s} (\psi_{\gamma,\delta}^{2})_{*}(\varphi_{\gamma,\delta}^{\prime}-\varphi_{\gamma,\delta})^{2} \mathrm{d}V
		:= J_{2,1} + J_{2,2}.
		\een
		By \eqref{lclfact1}, we have
		\ben \label{f-gamma-upper}
		(\varphi_{\gamma,\delta}^{2})_{*} = ((1-U^{\gamma/2+s}_{\delta})\mu^{\f12})^{2}_{*} \lesssim \delta^{2} \mu^{\f12}_{*}.
		\een
		Plugging \eqref{f-gamma-upper} into $J_{2,1}$ gives
		\ben \label{J21}
		J_{2,1}&\lesssim& \delta^{2}\delta^{-\gamma-2s}\int b_{s}(\theta)|v-v_{*}|^{-2s}  \mu^{\f12}_{*} (\psi_{\gamma,\delta}^{\prime}-\psi_{\gamma,\delta})^{2} \mathrm{d}V
		\\ \label{J21-line2}
		&=& \delta^{2}\delta^{-\gamma-2s} \mathcal{N}^{-2s,s}(\mu^{1/4},\psi_{\gamma,\delta})
		\lesssim \delta^{2}\delta^{-\gamma-2s}|U^{\gamma/2+s}_{\delta}(1-\rho \mu)f|^{2}_{\mathcal{L}^{s}_{-s}} \lesssim \delta^{2}|f|^{2}_{\mathcal{L}^{s}_{\gamma/2}},
		\een
		where 
		\beno
		\mathcal{N}^{\gamma,s}(g,h) := \int b_{s}(\theta)|v-v_{*}|^{\gamma}  g_{*}^{2} (h^{\prime}-h)^{2} \mathrm{d}V.
		\eeno
		In fact, the first inequality in \eqref{J21-line2} follows from the estimate (see \cite{he2022asymptotic})
		\beno
		\mathcal{N}^{\gamma,s}(\mu^{1/4},h) \lesssim |h|^{2}_{\mathcal{L}^{s}_{\gamma/2}},
		\eeno
		and the  second inequality follows from $|\mu|_{H^{2}} \lesssim 1$ and $\delta^{-\gamma/2-s}U^{\gamma/2+s}_{\delta} \in S^{\gamma/2+s}_{1,0}$ that is a radial  symbol of order $\gamma/2+s$.

		By  Lemma \ref{difference-term-complication-2}, we then have
		\ben \label{J22}
		J_{2,2}\lesssim\delta^{2s}\delta^{-\gamma-2s} \int (U^{\gamma/2+s}_{\delta}f)^{2}_{*}  \mathrm{d}v_{*}
		\lesssim \delta^{2s}|W_{\gamma/2+s}f|^{2}_{L^{2}}
		\lesssim \delta^{2s}|f|^{2}_{\mathcal{L}^{s}_{\gamma/2}}.
		\een
		Plugging the estimates \eqref{J21} and \eqref{J22} into \eqref{J2-part}, we get
		\ben \label{estimate-J-2}
		J_{2} \lesssim \delta^{2s}|f|^{2}_{\mathcal{L}^{s}_{\gamma/2}}.
		\een

		{\it \underline{Upper bound of $J_3$.}}
		Lemma \ref{difference-term-complication} gives
		\ben \label{estimate-J-3}
		J_{3} = \delta^{-\gamma-2s}\int b_{s}(\theta)|v-v_{*}|^{-2s} \left(h^{\prime} h^{\prime}_{*}-
		h h_{*}\right)^{2} \mu_{*}(1-\rho \mu)^{2}f^{2} \mathrm{d}V
		\lesssim  \delta^{2s} |f|^{2}_{L^{2}_{\gamma/2+s}}.
		\een

		{\it Step 3: The case when $-3s \leq \gamma <-2s$.}
		By plugging the estimates of $J_{1}$ in \eqref{estimate-J-1},  $J_{2}$ in \eqref{estimate-J-2}, $J_{3}$ in \eqref{estimate-J-3} into \eqref{separate-key-part},  for $0<\delta<1$, we have
		\beno   \mathcal{J}_{\rho,\gamma,s}(f) \gtrsim \delta^{-\gamma-2s}(1-\rho)^{\f72}|f|^{2}_{L^{2}_{\gamma/2+s}}-C\delta^{2s}|f|^{2}_{\mathcal{L}^{s}_{\gamma/2}}. \eeno
		Precisely, for some generic constants $0< C_{1} \leq 1 \leq  C_{2}$, we have
		\ben \label{key-estimate-order-s-2s}   \mathcal{J}_{\rho,\gamma,s}(f) \geq C_{1}(1-\rho)^{a}\delta^{s}|f|^{2}_{L^{2}_{\gamma/2+s}}-C_{2}\delta^{2s}|f|^{2}_{\mathcal{L}^{s}_{\gamma/2}}. \een
		In this step, $a=\f72$.
		By \eqref{J-rho-lower-bound} and \eqref{leading-with-l2-missing},  for some generic constants $0< C_{3} \leq 1 \leq  C_{4}$, we have
		\ben \label{known-estimate}    \mathcal{J}_{\rho,\gamma,s}(f) \geq C_{3} (1-\rho)^{2} |f|^{2}_{\mathcal{L}^{s}_{\gamma/2}}- C_{4}(1-\rho)^{\f12}|f|^{2}_{L^{2}_{\gamma/2+s}}.\een
		Then the  combination $\eqref{known-estimate} \times C_{5}\delta^{2s}(1-\rho)^{-2} + \eqref{key-estimate-order-s-2s}$ gives
		\ben \label{combination}  (1+C_{5}\delta^{2s}(1-\rho)^{-2}) \mathcal{J}_{\rho,\gamma,s}(f) &\geq& (
		C_{1}-C_{4}C_{5}(1-\rho)^{-\f32-a}\delta^{s})(1-\rho)^{a}\delta^{s}|f|^{2}_{L^{2}_{\gamma/2+s}}
		\\&&+(C_{3}C_{5}-C_{2})\delta^{2s}|f|^{2}_{\mathcal{L}^{s}_{\gamma/2}}. \nonumber \een
		We can then  take $C_{5}$ large enough such that $C_{3}C_{5}-C_{2} \geq C_{2}$, for example  $C_{5} = 2C_{2}/C_{3}$.
	 And  then we choose  $\delta$ small enough such that $C_{1}-C_{4}C_{5}(1-\rho)^{-\f32-a}\delta^{s} \geq 0$, for example  $\delta^{s}=\frac{C_{1}(1-\rho)^{\f32+a}}{C_{4}C_{5}} = \frac{C_{1}C_{3}(1-\rho)^{\f32+a}}{2C_{4}C_{2}}$. Thus,  we get
		\ben \label{conclusion1-gamma-minus2}   \mathcal{J}_{\rho,\gamma,s}(f) \geq C_{2} \delta^{2s}|f|^{2}_{\mathcal{L}^{s}_{\gamma/2}} = C_{2}\left(\frac{C_{1}C_{3}}{2C_{2}C_{4}}\right)^{2} (1-\rho)^{3+2a} |f|^{2}_{\mathcal{L}^{s}_{\gamma/2}}. \een
		Since $a=\f72$ in this step, we get $\mathcal{J}_{\rho,\gamma,s}(f) \gtrsim (1-\rho)^{10} |f|^{2}_{\mathcal{L}^{s}_{\gamma/2}}$
		for $-3s \leq \gamma < -2s$.

		{\it Step 4: The case when $-4s \leq \gamma < -3s$.}
		In this case, we take $\alpha = -3s, \beta=\gamma+3s$. Then $-s<\beta<0, \alpha+\beta=\gamma$. Similar to \eqref{separate-key-part},
		we get
		\ben \label{separate-gamma-very-small} \mathcal{J}_{\rho,\gamma,s}(f) &\gtrsim&
		\frac{1}{4}\delta^{-\beta}
		\int b_{s}(\theta)|v-v_{*}|^{\alpha} \mathbb{F}(\mu^{\f12}, U^{\beta/2}_{\delta}g) \mathrm{d}V
		\\  && - \frac{1}{2}\delta^{-\beta}\int b_{s}(\theta)|v-v_{*}|^{\alpha} \mathbb{F}((1-U^{\beta/2}_{\delta})\mu^{\f12}, U^{\beta/2}_{\delta}g) \mathrm{d}V
		\nonumber\\  && -  4\delta^{-\beta}\int b_{s}(\theta)|v-v_{*}|^{\alpha} \left(h^{\prime} h^{\prime}_{*}-
		h h_{*}\right)^{2} \mu_{*}g^{2} \mathrm{d}V
		:= \frac{1}{4} J^{\alpha,\beta}_{1}- \frac{1}{2} J^{\alpha,\beta}_{2} - 4 J^{\alpha,\beta}_{3},
		\nonumber
		\een
		where $h := U^{\beta/2}_{\delta}$.

		{\it \underline{ Lower bound of $J^{\alpha,\beta}_{1}$.}}
		Since $\alpha=-3s$, we can use the estimate \eqref{conclusion1-gamma-minus2} to obtain
		\ben \nonumber 
		 J^{\alpha,\beta}_{1} = \delta^{-\beta} \mathcal{J}_{\rho,-3s,s} (U^{\beta/2}_{\delta}f) 
		&\gtrsim& \delta^{-\beta}
		(1-\rho)^{10} |(\mathbb{I}-\mathbb{P}_{\rho})(U^{\beta/2}_{\delta}f)|^{2}_{\mathcal{L}^{s}_{-3s/2}} 
		\\ \nonumber  &\geq&
		\delta^{-\beta}
		(1-\rho)^{10} |(\mathbb{I}-\mathbb{P}_{\rho})(U^{\beta/2}_{\delta}f)|^{2}_{L^{2}_{-s/2}}
		\\ \label{J-al-be-1} &\gtrsim& \delta^{-\beta} (1-\rho)^{10} |f|^{2}_{L^{2}_{\gamma/2+s}},
		\een
		where we  have  used \eqref{delta-small-projection-to-not-projection} in the last inequality.

		{\it \underline{ Upper bound of $J^{\alpha,\beta}_{2}$.}} We now estimate 
		\beno J^{\alpha,\beta}_{2}= \delta^{-\beta}\int b_{s}(\theta)|v-v_{*}|^{\alpha} \mathbb{F}((1-U^{\beta/2}_{\delta})\mu^{\f12}, U^{\beta/2}_{\delta}(1-\rho \mu)f) \mathrm{d}V. \eeno
		For simplicity of notations,  set $\varphi_{\beta,\delta}= (1-U^{\beta/2}_{\delta})\mu^{\f12}, \psi_{\beta,\delta}=U^{\beta/2}_{\delta}(1-\rho \mu)f$. Then
		\ben \label{J-al-be-2-into-2}
		J^{\alpha,\beta}_{2} &\lesssim& \delta^{-\beta}\int b_{s}(\theta) |v-v_{*}|^{\alpha} (\varphi_{\beta,\delta}^{2})_{*}(\psi_{\beta,\delta}^{\prime}-\psi_{\beta,\delta})^{2} \mathrm{d}V \\ \nonumber &&+ \delta^{-\beta}\int b_{s}(\theta) |v-v_{*}|^{\alpha} (\psi_{\beta,\delta}^{2})_{*}(\varphi_{\beta,\delta}^{\prime}-\varphi_{\beta,\delta})^{2} \mathrm{d}V
		 \\ \nonumber &:=& J^{\alpha,\beta}_{2,1} + J^{\alpha,\beta}_{2,2}.
		\een
		Similar to \eqref{f-gamma-upper}, we get
		$ \varphi_{\beta,\delta}^{2}=((1-U^{\beta/2}_{\delta})\mu^{\f12})^{2} \lesssim \delta^{2} \mu^{\f12}. $
		This together with $\delta^{-\beta}U^{\beta/2}_{\delta} \leq W_{\beta/2}$ give
		\ben \label{J-al-be-2-into-2-1}
		J^{\alpha,\beta}_{2,1}&\lesssim& \delta^{2}\delta^{-\beta}\int b_{s}(\theta)|v-v_{*}|^{\alpha}  \mu^{\f12}_{*} (\psi_{\beta,\delta}^{\prime}-\psi_{\beta,\delta})^{2} \mathrm{d}V
		\nonumber \\&=& \delta^{2}\delta^{-\beta}\mathcal{N}^{-3s,s}(\mu^{1/4},\psi_{\beta,\delta}) \lesssim
		\delta^{2}\delta^{-\beta}|U^{\beta/2}_{\delta}f|^{2}_{\mathcal{L}^{s}_{-3s/2}} \lesssim \delta^{2}|f|^{2}_{\mathcal{L}^{s}_{\gamma/2}},
		\een
		which is similar to \eqref{J21-line2}.
		
		By Lemma \ref{difference-term-complication-2}, we have
		\ben \label{J-al-be-2-into-2-2}
		J^{\alpha,\beta}_{2,2} \lesssim
		\delta^{2s}\delta^{-\beta} \int (U^{\beta/2}_{\delta}f)^{2}_{*} \langle v_{*} \rangle^{\alpha+2s}  \mathrm{d}v_{*}
		\lesssim \delta^{2s}|W_{\gamma/2+s}f|^{2}_{L^{2}}
		\lesssim \delta^{2s}|f|^{2}_{\mathcal{L}^{s}_{\gamma/2}}.
		\een
		Plugging \eqref{J-al-be-2-into-2-1} and \eqref{J-al-be-2-into-2-2} into \eqref{J-al-be-2-into-2} gives
		\ben \label{J-al-be-2-into-2-final}
		J^{\alpha,\beta}_{2} \lesssim \delta^{2s}|f|^{2}_{\mathcal{L}^{s}_{\gamma/2}}.
		\een

		{\it \underline{ Upper bound of $J^{\alpha,\beta}_{3}$.}}
		By Lemma \ref{difference-term-complication}, we get
		\ben \label{J-al-be-3-up}
		J^{\alpha,\beta}_{3} = \delta^{-\beta}\int b_{s}(\theta)|v-v_{*}|^{\alpha} \left(h^{\prime} h^{\prime}_{*}-
		h h_{*}\right)^{2} \mu_{*}(1-\rho \mu)^{2}f^{2} \mathrm{d}V
		\lesssim  \delta^{2s} |f|^{2}_{L^{2}_{\gamma/2+s}}.
		\een

		Plugging the estimates of $J^{\alpha,\beta}_{1}$ in \eqref{J-al-be-1},
		$J^{\alpha,\beta}_{2}$ in \eqref{J-al-be-2-into-2-final},
		$J^{\alpha,\beta}_{3}$ in \eqref{J-al-be-3-up} into \eqref{separate-gamma-very-small},
		we get
		\beno   \mathcal{J}_{\rho,\gamma,s}(f) \gtrsim \delta^{-\beta} (1-\rho)^{10} |f|^{2}_{L^{2}_{\gamma/2+s}}-C\delta^{2s}|f|^{2}_{\mathcal{L}^{s}_{\gamma/2}}. \eeno
		Precisely, for some generic constants $C_{6}$ and $C_{7}$,
		we get
		\ben \label{key-estimate-by-He-2}   \mathcal{J}_{\rho,\gamma,s}(f) \geq C_{6}(1-\rho)^{10}\delta^{s}|f|^{2}_{L^{2}_{\gamma/2+s}}-C_{7}\delta^{2s}|f|^{2}_{\mathcal{L}^{s}_{\gamma/2}}. \een
		Together with coercivity estimate \eqref{known-estimate},
		by taking $a=10$ in \eqref{conclusion1-gamma-minus2}, for $-4s \leq \gamma < -3s$ it holds that 
		\ben \label{conclusion1-gamma-minus3}   \mathcal{J}_{\rho,\gamma,s}(f) \geq C_{7}\left(\frac{C_{6}C_{3}}{2C_{4}C_{7}}\right)^{2} (1-\rho)^{23} |f|^{2}_{\mathcal{L}^{s}_{\gamma/2}}. \een

		In general, we need $N_{\gamma,s}+1$ steps to reach $\gamma$ for $-3< \gamma<-2s$. By \eqref{conclusion1-gamma-minus2}, we define a sequence $\{a_n\}$ satisfying $a_{0}=\f72, a_{n}=3+2a_{n-1}$. This implies that  $a_{n}=13 \times 2^{n-1} -3$. Then for $-3< \gamma<-2s$,
		 by induction we will arrive at
		\ben \label{general-gamma-2-minus3}   \mathcal{J}_{\rho,\gamma,s}(f) \gtrsim (1-\rho)^{13 \times 2^{N_{\gamma,s}} -3} |f|^{2}_{\mathcal{L}^{s}_{\gamma/2}}. \een
	This completes the proof of the theorem.
	\end{proof}

	{\bf Acknowledgments.} 
	The research was partially supported by the National Key Research and
	Development Program of China project no. 2021YFA1002100. The research of Tong Yang was
	supported by a fellowship award from the Research Grants Council of the Hong Kong Special
	Administrative Region, China (Project no. SRF2021-1S01). The research of Yu-Long Zhou was
	supported by the NSFC project no. 12001552, the Science and Technology Project in Guangzhou
	no. 202201011144.
	
	\bibliographystyle{siam}
	\bibliography{quantum-B-IPL}

\begin{thebibliography}{10}

\bibitem{alexandre2019global}
{\sc R.~Alexandre, F.~H{\'e}rau, and W.-X. Li}, {\em {Global hypoelliptic and
  symbolic estimates for the linearized Boltzmann operator without angular
  cutoff}}, Journal de Math{\'e}matiques Pures et Appliqu{\'e}es, 126 (2019),
  pp.~1--71.

\bibitem{alexandre2011boltzmann}
{\sc R.~Alexandre, Y.~Morimoto, S.~Ukai, C.-J. Xu, and T.~Yang}, {\em {The
  Boltzmann equation without angular cutoff in the whole space: II, Global
  existence for hard potential}}, Analysis and Applications, 9 (2011),
  pp.~113--134.

\bibitem{alexandre2012boltzmann}
\leavevmode\vrule height 2pt depth -1.6pt width 23pt, {\em {The Boltzmann
  equation without angular cutoff in the whole space: I, Global existence for
  soft potential}}, Journal of Functional Analysis, 262 (2012), pp.~915--1010.

\bibitem{bae2021relativistic}
{\sc G.-C. Bae, J.~W. Jang, and S.-B. Yun}, {\em {The relativistic quantum
  Boltzmann equation near equilibrium}}, Archive for Rational Mechanics and
  Analysis, 240 (2021), pp.~1593--1644.

\bibitem{baranger2005explicit}
{\sc C.~Baranger and C.~Mouhot}, {\em {Explicit spectral gap estimates for the
  linearized Boltzmann and Landau operators with hard potentials}}, Revista
  Matem{\'a}tica Iberoamericana, 21 (2005), pp.~819--841.

\bibitem{benedetto2004some}
{\sc D.~Benedetto, F.~Castella, R.~Esposito, and M.~Pulvirenti}, {\em {Some
  considerations on the derivation of the nonlinear quantum Boltzmann
  equation}}, Journal of Statistical Physics, 116 (2004), pp.~381--410.

\bibitem{benedetto2006some}
\leavevmode\vrule height 2pt depth -1.6pt width 23pt, {\em {Some considerations
  on the derivation of the nonlinear quantum Boltzmann equation II: the low
  density regime}}, Journal of Statistical Physics, 124 (2006), pp.~951--996.

\bibitem{benedetto2007short}
\leavevmode\vrule height 2pt depth -1.6pt width 23pt, {\em {A short review on
  the derivation of the nonlinear quantum Boltzmann equations}}, Communications
  in Mathematical Sciences, 5 (2007), pp.~55--71.

\bibitem{benedetto2008n}
\leavevmode\vrule height 2pt depth -1.6pt width 23pt, {\em {From the N-body
  Schr{\"o}dinger equation to the quantum Boltzmann equation: a term-by-term
  convergence result in the weak coupling regime}}, Communications in
  Mathematical Physics, 277 (2008), pp.~1--44.

\bibitem{benedetto2005weak}
{\sc D.~Benedetto, M.~Pulvirenti, F.~Castella, and R.~Esposito}, {\em {On the
  weak-coupling limit for bosons and fermions}}, Mathematical Models and
  Methods in Applied Sciences, 15 (2005), pp.~1811--1843.

\bibitem{briant2016cauchy}
{\sc M.~Briant and A.~Einav}, {\em {On the Cauchy problem for the homogeneous
  Boltzmann--Nordheim equation for bosons: local existence, uniqueness and
  creation of moments}}, Journal of Statistical Physics, 163 (2016),
  pp.~1108--1156.

\bibitem{chapman1990mathematical}
{\sc S.~Chapman and T.~G. Cowling}, {\em {The mathematical theory of
  non-uniform gases: an account of the kinetic theory of viscosity, thermal
  conduction and diffusion in gases}}, Cambridge university press, 1990.

\bibitem{chen2021emergence}
{\sc T.~Chen and M.~Hott}, {\em {On the emergence of quantum Boltzmann
  fluctuation dynamics near a Bose-Einstein Condensate}}, arXiv preprint
  arXiv:2106.03825,  (2021).

\bibitem{duan2022solutions}
{\sc R.~Duan, L.-B. He, T.~Yang, and Y.-L. Zhou}, {\em Solutions to the
  non-cutoff boltzmann equation in the grazing limit}, Annales de l'Institut
  Henri Poincar{\'e} C,  (2022).

\bibitem{erdHos2004quantum}
{\sc L.~Erd{\H{o}}s, M.~Salmhofer, and H.-T. Yau}, {\em {On the quantum
  Boltzmann equation}}, Journal of Statistical Physics, 116 (2004),
  pp.~367--380.

\bibitem{escobedo2007fundamental}
{\sc M.~Escobedo, S.~Mischler, and J.~J.~L. Vel{\'a}zquez}, {\em {On the
  fundamental solution of a linearized Uehling--Uhlenbeck equation}}, Archive
  for Rational Mechanics and Analysis, 186 (2007), pp.~309--349.

\bibitem{escobedo2008singular}
\leavevmode\vrule height 2pt depth -1.6pt width 23pt, {\em {Singular solutions
  for the Uehling--Uhlenbeck equation}}, Proceedings of the Royal Society of
  Edinburgh Section A: Mathematics, 138 (2008), pp.~67--107.

\bibitem{gressman2011global}
{\sc P.~Gressman and R.~M. Strain}, {\em {Global classical solutions of the
  Boltzmann equation without angular cut-off}}, Journal of the American
  Mathematical Society, 24 (2011), pp.~771--847.

\bibitem{guo2003classical}
{\sc Y.~Guo}, {\em {Classical solutions to the Boltzmann equation for molecules
  with an angular cutoff}}, Archive for Rational Mechanics and Analysis, 169
  (2003), pp.~305--353.

\bibitem{he2021boltzmann}
{\sc L.-B. He and Y.-L. Zhou}, {\em {Boltzmann equation with cutoff Rutherford
  scattering cross section near Maxwellian}}, Archive for Rational Mechanics
  and Analysis, 242 (2021), pp.~1631--1748.

\bibitem{he2022asymptotic}
\leavevmode\vrule height 2pt depth -1.6pt width 23pt, {\em {Asymptotic analysis
  of the linearized Boltzmann collision operator from angular cutoff to
  non-cutoff}}, Annales de l'Institut Henri Poincar{\'e} C, 39 (2022),
  pp.~1097--1178.

\bibitem{lerner2013phase}
{\sc N.~Lerner, Y.~Morimoto, K.~Pravda-Starov, and C.-J. Xu}, {\em {Phase space
  analysis and functional calculus for the linearized Landau and Boltzmann
  operators}}, Kinetic \& Related Models, 6 (2013), pp.~625--648.

\bibitem{li2019global}
{\sc W.~Li and X.~Lu}, {\em {Global existence of solutions of the Boltzmann
  equation for Bose--Einstein particles with anisotropic initial data}},
  Journal of Functional Analysis, 276 (2019), pp.~231--283.

\bibitem{lu2000modified}
{\sc X.~Lu}, {\em {A modified Boltzmann equation for Bose--Einstein particles:
  isotropic solutions and long-time behavior}}, Journal of Statistical Physics,
  98 (2000), pp.~1335--1394.

\bibitem{lu2004isotropic}
\leavevmode\vrule height 2pt depth -1.6pt width 23pt, {\em {On isotropic
  distributional solutions to the Boltzmann equation for Bose-Einstein
  particles}}, Journal of Statistical Physics, 116 (2004), pp.~1597--1649.

\bibitem{lu2005boltzmann}
\leavevmode\vrule height 2pt depth -1.6pt width 23pt, {\em {The Boltzmann
  equation for Bose--Einstein particles: velocity concentration and convergence
  to equilibrium}}, Journal of Statistical Physics, 119 (2005), pp.~1027--1067.

\bibitem{lukkarinen2009not}
{\sc J.~Lukkarinen and H.~Spohn}, {\em {Not to normal order--notes on the
  kinetic limit for weakly interacting quantum fluids}}, Journal of Statistical
  Physics, 134 (2009), pp.~1133--1172.

\bibitem{mouhot2006explicit}
{\sc C.~Mouhot}, {\em {Explicit coercivity estimates for the linearized
  Boltzmann and Landau operators}}, Communications in Partial Differential
  Equations, 31 (2006), pp.~1321--1348.

\bibitem{mouhot2007spectral}
{\sc C.~Mouhot and R.~M. Strain}, {\em {Spectral gap and coercivity estimates
  for linearized Boltzmann collision operators without angular cutoff}},
  Journal de Math{\'e}matiques Pures et Appliqu{\'e}es, 87 (2007),
  pp.~515--535.

\bibitem{nordhiem1928kinetic}
{\sc L.~Nordhiem}, {\em {On the kinetic method in the new statistics and
  application in the electron theory of conductivity}}, Proceedings of the
  Royal Society of London. Series A, Containing Papers of a Mathematical and
  Physical Character, 119 (1928), pp.~689--698.

\bibitem{ouyang2022quantum}
{\sc Z.~Ouyang and L.~Wu}, {\em On the quantum boltzmann equation near
  maxwellian and vacuum}, Journal of Differential Equations, 316 (2022),
  pp.~471--551.

\bibitem{pao1974boltzmann1}
{\sc Y.-P. Pao}, {\em {Boltzmann collision operator with inverse-power
  intermolecular potentials, I}}, Communications on Pure and Applied
  Mathematics, 27 (1974), pp.~407--428.

\bibitem{pao1974boltzmann2}
\leavevmode\vrule height 2pt depth -1.6pt width 23pt, {\em {Boltzmann collision
  operator with inverse-power intermolecular potentials, II}}, Communications
  on Pure and Applied Mathematics, 27 (1974), pp.~559--581.

\bibitem{uehling1933transport}
{\sc E.~A. Uehling and G.~Uhlenbeck}, {\em {Transport phenomena in
  Einstein-Bose and Fermi-Dirac gases. i}}, Physical Review, 43 (1933), p.~552.

\bibitem{wang1952propagation}
{\sc C.~Wang~Chang and G.~Uhlenbeck}, {\em On the propagation of sound in
  monatomic gases}, Studies in Statistical Mechanics,  (1952), pp.~43--75.

\bibitem{zhou2022global}
{\sc Y.-L. Zhou}, {\em {Global well-posedness of the quantum Boltzmann equation
  for bosons interacting via inverse power law potentials}}, arXiv preprint
  arXiv:2210.08428,  (2022).

\end{thebibliography}

\end{document}